\documentclass[conference]{IEEEtran}
\IEEEoverridecommandlockouts
\usepackage{cite}
\usepackage{amsmath,amssymb,amsfonts,color,graphicx,amsthm,multicol,multirow,booktabs,url,lastpage}
\usepackage{pstricks,pst-sigsys,pst-circ,pst-plot,pst-eucl,pstricks-add}
\usepackage{float}
\usepackage{algorithmic}
\usepackage{graphicx}
\usepackage{textcomp}
\usepackage{xcolor}
\def\BibTeX{{\rm B\kern-.05em{\sc i\kern-.025em b}\kern-.08em
    T\kern-.1667em\lower.7ex\hbox{E}\kern-.125emX}}

\usepackage{caption}
\usepackage{subcaption}

\usepackage{cleveref}
\crefname{claim}{Claim}{claims}
\crefname{theorem}{Theorem}{theorems}
\crefname{definition}{Definition}{definitions}
\crefname{lemma}{Lemma}{lemmas}
\crefname{fact}{Fact}{facts}
\crefname{corollary}{Corollary}{corollaries}
\crefname{coroc}{Remark}{remarks}
\crefname{corot}{Remark}{remarks}
\crefname{corod}{Remark}{remarks}
\crefname{corol}{Remark}{remarks}
\crefname{figure}{Figure}{fig}

\newcommand{\indep}{\rotatebox[origin=c]{90}{$\models$}}

\newcommand{\sgn}[1]{\operatorname*{sgn}\left( {#1} \right)}
\newcommand{\tr}[1]{\operatorname*{tr}\left( {#1} \right)}
\newcommand{\bs}[1]{\boldsymbol{#1}}
\newcommand{\prob}[1]{\mathbb{P}\left\lbrace {#1} \right\rbrace}
\newcommand{\expc}[3]{\mathbb{E}_{#1}^{#2} \left[ {#3} \right]}
\newcommand{\var}[3]{\mathbb{V}ar_{#1}^{#2} \left[ {#3} \right]}
\newcommand{\dive}[3]{D_{#1} \left( {#2} || {#3} \right)}
\newcommand{\DR}[2]{%
  \mathrel{\mathop\gtrless\limits^{#1}_{#2}}%
}

\newtheorem{definition}{Definition}
\newtheorem{theorem}{Theorem}
\newtheorem{claim}{Claim}
\newtheorem{examplee}{Example}

\newtheorem{lemma}{Lemma}
\newtheorem{corollary}{Corollary}[claim]
\newtheorem{coroc}{Remark}[claim]

\newtheorem{corol}{Remark}[lemma]

\begin{document}

\title{Covertly Controlling a Linear System\\
}
\author{\IEEEauthorblockN{Barak Amihood}
\IEEEauthorblockA{\textit{The School of Electrical and Computer Engineering} \\
\textit{Ben-Gurion University of the Negev}\\
Beer-Sheva, Israel \\
barakam@post.bgu.ac.il}
\and
\IEEEauthorblockN{Asaf Cohen}
\IEEEauthorblockA{\textit{The School of Electrical and Computer Engineering} \\
\textit{Ben-Gurion University of the Negev}\\
Beer-Sheva, Israel \\
coasaf@bgu.ac.il}
}

\maketitle

\begin{abstract}
Consider the problem of covertly controlling a linear system. In this problem, Alice desires to control (stabilize or change the parameters of) a linear system, while keeping an observer, Willie, unable to decide if the system is indeed being controlled or not. 

We formally define the problem, under two different models: (i) When Willie can only observe the system's output (ii) When Willie can directly observe the control signal. Focusing on AR(1) systems, we show that when Willie observes the system's output through a clean channel, an inherently unstable linear system can not be covertly stabilized. However, an inherently stable linear system can be covertly controlled, in the sense of covertly changing its parameter. Moreover, we give direct and converse results for two important controllers: a minimal-information controller, where Alice is allowed to used only $1$ bit per sample, and a maximal-information controller, where Alice is allowed to view the real-valued output. Unlike covert communication, where the trade-off is between rate and covertness, the results reveal an interesting \emph{three--fold} trade--off in covert control: the amount of information used by the controller, control performance and covertness. To the best of our knowledge, this is the first study formally defining covert control.
\end{abstract}


\section{Introduction}

The main objective in control theory is to develop algorithms that govern or control systems. Usually, the purpose of a selected algorithm is to drive the system to a desired state (or keeping it at a certain range of states) while adhering to some constraints such as rate of information, delay, power and overshoot. Essentially, ensuring a level of control (under some formal definition) subject to some (formally defined) constraint. Traditionally, the controller monitors the system's process and compares it with a reference. The difference between the actual and desired value of the process (the error signal) is analysed and processed in order to generate a control action, which in tern is applied as a feedback to the system, bringing the controlled process to its desired state.

Numerous systems require external control in order to operate smoothly. For example, various sensors in our power, gas or water networks monitor the flow and control it. We use various signaling methods to control cameras in homeland security applications or medical devices. In some cases, the controlling signal is manual, while in others it is an automatic signal inserted by a specific part of the application in charge of control. In this paper, however, we extend the setup to cases where achieving control over the system is not the only goal, and consider the control problem while adding an additional constraint of covertness: either staying undetected by an observer (taking the viewpoint of an illegitimate controller), or being able to detect any control operation (taking the viewpoint of a legitimate system owner). This additional covertness constraint clearly seems reasonable in applications such as security or surveillance systems, but, in fact, recent attacks on medical devices \cite{mahler2018know} and civil infrastructures \cite{doi:10.1177/0361198118756885} stress out the need to account for covert control in a far wider spectrum of applications. As a result, it is natural to ask: Can one covertly control a linear system? If so, in what sense? When is it possible to identify with high probability any attempt to covertly control a system? Do answers to the above question depend on the type of controller used, in terms of complexity or information gathered from the system?

While the problem we define herein is related to covert communication, which was studied extensively in the information theory community \cite{bash_limits,Reliable_Deniable,Achieving_Undetectable_Communication,resolvability_perspective,First_and_Second_Order_Asymptotics_in_Covert_Communication,9035417}, key differences immediately arise. First, in covert communication, Alice desires to covertly send a message to Bob. Hence, fixing a covertness constraint, Alice's success is measured in rate - the number of bits per channel use she is able to send covertly. Herein, Alice's objective is \emph{control}, which we measure in the ability to stabilize the system, or change its parameters. Second, in the covert control setup, we identify an additional dimension, which does not exist in covert communication: the amount of information Alice has to extract from the system in order to carry out her objective (covertly control the system). Thus, it is clear there is a non-trivial interplay between information, control and covertness, adding depth and a multitude of open problems. 
\subsection{Main Contribution}
\begin{figure}
\centering
\includegraphics[width=0.8\linewidth]{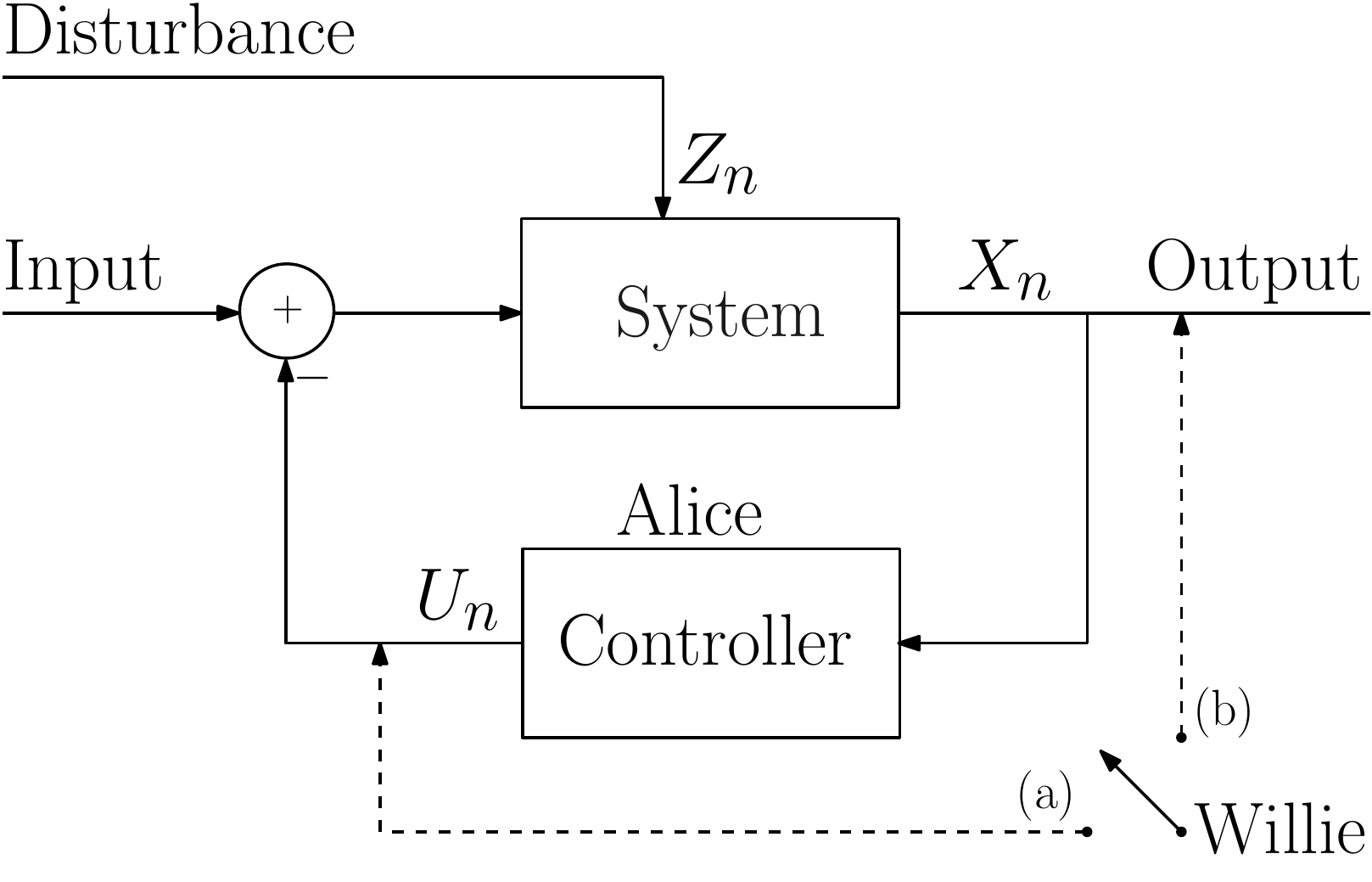}
\caption{A basic Covert Control model. 
Alice observes the system's output, and wishes to control the system, either using a desired reference signal or without it. She can decide how frequently and how accurately to sample the output on the one hand, and how frequent and how strong will be her control signal be, on the other. The control signal is an input to the system, but it is also observed by Willie either explicitly (a), or implicitly when observing the output of the system (b). In both cases, Willie's observations can be either clean, or viewed via a noisy channel. Alice's goal is to control the system without being noticed. Willie's goal is to detect if Alice is indeed controlling the system.}
\label{fig:Covert_Control}
\end{figure}

In this paper, we focus on a simple linear system, schematically depicted in \Cref{fig:Covert_Control}. Alice, observing the system's output $X_n$, wishes to control the system's behaviour though a control signal $U_n$. We assume the system without control follows a first order Auto Regressive model (AR(1)), hence,
\begin{align}\label{eq:simple_linear_model}
\begin{split}
X_{n+1} = a X_{n} + Z_{n} - U_{n},
\end{split}
\end{align}
and focus on two types of control objectives: (i) Stabilizing an otherwise unstable system. (ii) Changing the parameter ($a$) of a stable system. However, Alice wishes to perform her control action without being detected by Willie. We discuss two different scenarios in terms of Willie's observations. (i) Willie is able to observe Alice's control signal. (ii) Willie can only observe the system's output. In both cases, Willie's observations might be clean, or through a noisy channel. 

We formally define covertness and detection in these scenarios. We then show that an unstable AR(1) cannot be covertly controlled, in the sense of maintaining a finite $\gamma$-moment without being detected. On the other hand, we show that a Gaussian stable AR(1) system can be covertly manipulated by Alice, in the sense of changing its gain without being detected by Willie. We then turn to specific results for two interesting, limiting--cases controllers. The first is a \emph{minimum--information} controller, in which Alice is allowed to retrieve only one bit per sample. We show that under this strict information constraint, Alice cannot covertly control the system, no matter if Willie observes the control signal directly or only the system's output. On the other hand, we then turn to the opposite limiting case, a \emph{maximum--information} controller, in which Alice is allowed to view the real--valued signal. We show that with maximum information, if the system's gain is small enough, Alice can manipulate a system, ``resetting" it to its initial value at a time instance of her choice, without being detected by Willie. Conversely, for a large enough system gain, Alice cannot reset the system without being detected by Willie.

The results above, to the best of our knowledge, are the first to characterize cases where one can or cannot covertly control a linear system. They reveal an important three--fold trade--off between information, covertness and control, yet to be fully characterized.

\subsection{Related Work}
\label{subsection:Scientific_overview}
\subsubsection{Covert Communication}Communicating covertly has been a long standing problem. In this scenario, two parties, Alice and Bob, wish to communicate, while preventing a third party, Willie, from detecting the mere presence
of communication. While studied in the steganography and spread spectrum communication literature \cite{Steganography, Spread, bash_hiding}, the first information theoretic investigation was done in \cite{bash_limits, Square_Root_Law}. This seminal work considered Additive White Gaussian Noise channels (whose variances are known to all parties), and computed the highest achievable rate between Alice and Bob, while ensuring Willie's sum of false alarm and missed detection probabilities ($\alpha + \beta$) is arbitrarily close to one. 
This, however, resulted in a transmission rate which is asymptotically negligible, e.g., $O(\sqrt{n})$ bits for $n$ channel uses. In fact, this ``square root law" for covert communication holds more generally, e.g., for Binary Symmetric Channels \cite{Reliable_Deniable} and via channel resolvability \cite{resolvability_perspective}. These works strengthened the insight that forcing $\alpha + \beta$ arbitrarily close to $1$, and assuming Willie uses optimal detection strategies, results in a vanishing rate. In order to achieve a strictly positive covert communication rate, Alice and Bob need some advantage over Willie
\cite{Not_Know_Noise, Not_Know_When, Reliable_deniable_communication_with_channel_uncertainty, Achieving_Undetectable_Communication, Achieving_positive_rate_with_undetectable_communication_over_AWGN_and_Rayleigh_channels, Achieving_positive_rate_with_undetectable_communication_Over_MIMO_rayleigh_channels}, or, alternatively, intelligently use the fact that Willie cannot use an optimal detector in practice, or some key problem parameters are beyond his reach.
For example, \cite{Not_Know_Noise} showed that if Willie does not know his exact noise statistics, Alice can covertly transmit $O(n)$ bits over $T$ transmission slots, one of which she can utilize, where each slot duration can encompass $n$ bit codewords. 
The authors in \cite{Not_Know_When} showed that if Alice and Bob secretly pre-arrange on which out of $T(n)$ slots Alice is going to transmit, they can reliably exchange $O(\sqrt{n \log{T(n)}})$ bits on an AWGN channel covertly from Willie.

A critical aspect of understanding communication systems under practical constraints is their analysis and testing in finite block length regimes. While asymptotic behaviour gives us fundamental limits and important insights, it is critical to understand such systems with finite, realistic block lengths, demanded by either complexity or delay constraints. The first studies were in 
\cite{Delay_Intolerant_Covert_Communications_With_Either_Fixed_or_Random_Transmit_Power, Covert_communication_with_finite_blocklength_in_AWGN_channels, Covert_communications_with_extremely_low_power_under_finite_block_length_over_slow_fading, Delay_Constrained_Covert_Communications_With_a_Full_Duplex_Receiver, A_Finite_Block_Length_Achievability_Bound_for_Low_Probability_of_Detection_Communication, First_and_Second_Order_Asymptotics_in_Covert_Communication, Finite_Blocklength_Analysis_of_Gaussian_Random_coding_in_AWGN_Channels_under_Covert_constraints_II_Viewpoints_of_Total_Variation_Distance}. 
Needless to say, practical constrains such as limited delay and finite blocks, may limit the communication of Alice and Bob, yet, on the other hand, may limit the warden's ability to detect the communication, hence it is not a priori clear which will have a stronger effect.

Several studies considered a model in which Alice and Bob can utilize a friendly helper which
can generate Artificial Noise (AN) 
\cite{Covert_Communication_in_the_Presence_of_an_Uninformed_Jammer, Covert_Wireless_Communication_With_Artificial_Noise_Generation, Covert_Communications_with_a_Full_Duplex_Receiver_over_Wireless_Fading_Channels, Passive_Self_Interference_Suppression_for_Full_Duplex_Infrastructure_Nodes}. 
In those models, the jammer is consider either to be an ``outsider", or just Bob utilizing an additional antenna for
transmitting AN. 
Those studies have shown an improvement of the covert communication rate between Alice and Bob compared to cases without a jammer.

\subsubsection{Control}In many contemporary applications, however, communication is used merely as means to accomplish a certain task, namely, communication is needed since the data required for the task is located on various (remote) sensors, devices or network locations. The utility in such applications depends on the task and thus, the tension when adding covertness constraints is not necessarily just between covertness and the number of bits transmitted, the task itself must be accounted for. For example, in the context of covertness when controlling a system, the tension is not just between the rate of information gathered from the system and covertness -- other aspects come into play, such as stability, delay, etc.

Focusing on control applications, it is not a-priori clear where to measure the control signal, and what are the interesting trade--offs to characterize. For example, in \cite{kostina_exact} the authors considered the minimum number of bits to stabilize a linear system, i.e., the system in \cref{eq:simple_linear_model}, where $\{Z_{n}\}$ were assumed independent random variables, with bounded $\alpha$-th moments, and  $\{U_{n}\}$ were the control actions, chosen by a controller, who received each time instant a single element of a finite set $\{ 1, \ldots, M \}$, as its only information about system state. The authors showed that for $|a| > 1$ (an inherently unstable system), $M = \lfloor a \rfloor + 1$ is necessary and sufficient to achieve $\beta$-moment stability, for any $\beta < \alpha$. 
Their approach was to use a normal/emergency (zoom-in/zoom-out) based controller, in order to manage the magnitude of $X_{n}$. This control signal, however, was not meant to be undetectable, and is indeed far from covert. Similar to \cite{kostina_exact},\cite{Stabilizing_a_system_with_an_unbounded_random_gain, rate-cost, Control_Over_Gaussian, Minimum_Data_Rate} 
also considered controlling a linear system, under different constraints, such as rate and energy. Yet again, covertness was not taken into account in these works.

\section{Problem formulation}\label{section:Problem_formulation}
In the context of \cref{fig:Covert_Control} and \cref{eq:simple_linear_model},  this work is on the following two cases. 
(i) Willie is observing Alice's controller's output. (ii) Willie is observing Alice's system's output. 
In each of the cases, Willie tries to detect any controlling action on Alice's part, by observing her controller's output or system's output, through some channel (clean or noisy). In turn, Alice's will is to control her system, while staying undetected by Willie. 
Alice's system is modeled as an AR(1) system (\cref{eq:simple_linear_model}), where $X_{0} = 0$, $a$ is the system \emph{gain}, $\{ Z_{n} \}$ are identically distributed, independent random variables, and $U_{n}$ is the control action at time $n$.
Choosing this kind of model stems from the fact that an AR(1) process is well-known, well studied, and can constitute a simple model, on the one hand, yet can capture the complexity of the problem on the other. It can model systems with memory, having both random input,  and a built-in system's gain. 

The control signal, which we denote by $U_{n}$, can function in several manners. First, Alice can use it to stabilize the system, with respect to some stabilization criteria. Second, Alice can use it to alter the system's parameters or to intervene with the system's operation process.
For simplicity, we focus on controllers of the form $U_{n} = f(X_{n-1})$.

Any controller of the form $U_{n} = f(X_{n-1})$ can use different amount of information from $X_{n-1}$ in order to control the system. This amount is determined by the function $f(\cdot)$. 
Explicitly, $f(\cdot)$ can use all the information in $X_{n-1}$, i.e., use infinite number of bits representation of $X_{n-1}$, or use some quantized value of $X_{n-1}$, i.e., use a finite number of bits in order to represent $X_{n-1}$. In this work, besides the general results, we also consider two different specific controllers, which give a glimpse on the following two extreme cases. First, the case in which a minimal amount of information is taken from $X_{n-1}$, i.e., a single bit per sample. 
Second, the case in which maximal information is taken from $X_{n-1}$, i.e., an infinite number of bits. 
Those two controllers are given in 
\cref{definition:One_bit_controller} and 
\cref{definition:Threshold_controller} below, respectively.

For an AR(1) system with limited support noises, i.e., $\{ Z_{n }\}$ for which there is a $B > 0$, such that for any $n$, $|Z_{n}| \leq B$, we define the following controller.
\begin{definition}\label{definition:One_bit_controller}
Let $U_{n}^{one}$ be the following control signal \cite{kostina_exact} 
\begin{align}\label{eq:One_bit_controller}
U_{n}^{one} = \frac{a}{2} C_{n-1} \operatorname*{sgn}{(X_{n-1})},
\end{align}
where $B$ is the noise bound and $a \in \mathbb{R}$, $a \neq 2$. $C_n$ is given by the following recursive formula, $C_{n} = (a/2) C_{n-1} + B$, and $C_{1} \geq \frac{B}{1- a/2}$. We refer to $U_{n}^{one}$ as the one-bit controller.
\end{definition}

Note that this controller keeps $X_{n}$ bounded, but only for limited support noise \cite{kostina_exact}.
\begin{definition}\label{definition:Threshold_controller}
Let $U_{n}^{th}$ be the following control signal
\begin{align}\label{eq:Threshold_controller}
U_{n}^{th} = 
\left\lbrace
\begin{matrix}
a X_{n-1}, & |X_{n-1}| \geq D \\
0, & |X_{n-1}| < D 
\end{matrix}
\right. ,
\end{align}

where $D > 0$ is a threshold value and $a$ is the system's gain. We refer to $U_{n}^{th}$ as the threshold controller.
\end{definition}

This intuitive controller acts as a reset to the system, that is, when the system state crosses some level, which we denote as $D$, the system returns to the initial state $X_{0} = 0$. 

In order to measure the covertness of Alice's control policy, and the detectability at Willie side, we define a covertness criterion and a detection criterion. 
Specifically, we introduce the notion of $\varepsilon$-covertness and $1-\delta$-detection.
Denote by $\alpha$ and $\beta$, the false alarm and  miss detection probabilities for hypothesis testing problem, respectively.
\begin{definition}\label{definition:Epsilon_covertness}
We say that $\varepsilon$-covertness is achieved by Alice, if for some $\varepsilon > 0$ we have, $\alpha + \beta \geq 1 - \varepsilon$.
\end{definition}

\cref{definition:Epsilon_covertness} is a well-known criteria in covert communication, first used in  \cite{bash_limits} to establish the fundamental limit
of covert communication, and next used ubiquitous, e.g., 
\cite{Reliable_Deniable,Not_Know_Noise,Not_Know_When,dvorkind_maximizing,dvorkind_rate}. 

\begin{definition}\label{definition:1-delta_detection}
We say that $1-\delta$-detection is achieved by Willie, if for some $\delta > 0$ we have, $\alpha + \beta \leq \delta$.
\end{definition}



\subsection{Preliminaries}
\label{subsection:Preliminary}
In this sub-section, we introduce some relevant preliminary results and additional definitions.



We call an AR(1) process with $\{ Z_{n} \}$ being a white Gaussian noise (WGN), i.e., $\{ Z_{n} \} \sim \mathcal{N}(0, \sigma_{Z}^{2})$, a \emph{Gaussian AR(1) process}.
\begin{lemma}\label{lemma:PDF_of_n_samples_AR_1_process}
Let $\textbf{X}^{(n)} \triangleq 
\begin{pmatrix}
X_{1} & X_{2} & \cdots & X_{n}
\end{pmatrix}$ 
be an $n$-tupple of a Gaussian AR(1) process. Then, the PDF of 
$\textbf{X}^{(n)}$ is given by,
\begin{align}\label{eq:PDF_of_n_samples_AR_1_process}
\begin{split}
f_{\textbf{X}^{(n)}}(\textbf{x})
=
\frac{1}{\left( 2 \pi \right)^{\frac{n}{2}} \sqrt{|\bs{\Sigma}|}} e^{-\frac{1}{2} \textbf{x}^{T} \bs{\Sigma}^{-1} \textbf{x}},
\end{split}
\end{align}

where 
$\bs{\Sigma} = \textbf{A} \textbf{A}^{T}$ is a square full-rank matrix, and $\textbf{A}$ is defined as $\textbf{A}_{i,j} = \sigma_{Z} a^{i-j} u(i-j)$, where $u(\cdot)$ is the discrete step function.
I.e., $\textbf{X}^{(n)} \sim \mathcal{N}(\textbf{0}, \bs{\Sigma})$. In addition, $[\bs{\Sigma}]_{ij} = \frac{\sigma_{Z}^{2}}{1 - a^{2}} \left( a^{|i-j|} - a^{i+j} \right), \quad \forall 1 \leq i,j \leq n$. 
\end{lemma}

\begin{corol}\label{corol:Stable_AR(1)_Covariance_in_SS}
Under steady-state, and for $|a| < 1$, we have, 
\[
\bs{\Sigma}_{ij}
=
\frac{\sigma_{Z}^{2}}{1 - a^{2}} a^{|i-j|}, \quad \forall 1 \leq i,j \leq n,
\]

which yields a wide-sense stationary process and a stable system. 
\end{corol}

\begin{definition}[\cite{THOMAS_ch2}]
For two distributions defined over the same support, $P$ and $Q$, the Kullback-Leibler Divergence is defined by
\begin{align}\label{eq:Kullback-Leibler_divergence}
\begin{split}
\dive{KL}{P}{Q}
\triangleq
\expc{P}{}{\log{\left( \frac{P(x)}{Q(x)} \right)}}.
\end{split}
\end{align}
\end{definition}

\begin{definition}[\cite{THOMAS_ch11}] 
For two distributions defined over the same support, $P$ and $Q$,
the total variation is defined by
\begin{align}\label{eq:Total_variation}
\begin{split}
\mathcal{V}_{T}(P,Q)
\triangleq
\frac{1}{2} \expc{Q}{}{\left\vert \frac{P(x)}{Q(x)} - 1 \right\vert}.
\end{split}
\end{align}
\end{definition}

\begin{lemma}[\cite{THOMAS_ch11}]\label{lemma:Connection_between_divergence_to_total_variation} For any two distributions defined over the same support, $P$ and $Q$, we have 
\begin{align}\label{eq:Connection_between_divergence_to_total_variation}
\begin{split}
\mathcal{V}_{T}(P,Q)
\leq
\sqrt{\frac{1}{2} \dive{KL}{P}{Q}}.
\end{split}
\end{align}
\end{lemma}
\begin{lemma}[\cite{Lehmann2005}]\label{lemma:Connection_between_error_probabilities_to_total_variation}
 For the optimal test, the sum of the error probabilities, i.e., $\alpha + \beta$, is given by
\begin{align}\label{eq:Connection_between_error_probabilities_to_total_variation}
\begin{split}
\alpha + \beta
=
1 - \mathcal{V}_{T}(P,Q).
\end{split}
\end{align}
\end{lemma}
\begin{lemma}\label{lemma:Divergence_of_gaussian_vector}
Assume two multivariate normal distributions, $\mathcal{N}_{0}$ and $\mathcal{N}_{1}$, with means $\bs{\mu} _{0}$ and $\bs{\mu} _{1}$, of the same dimension $n$, and with (non-singular) covariance matrices, $\bs{\Sigma}_{0}$ and $\bs{\Sigma}_{1}$, respectively. 
Then, the Kullback–Leibler divergence between the distributions is 
\begin{multline}\label{eq:Divergence_of_gaussian_vector}
D_{KL} \left( 
\mathcal{N}_{0} || \mathcal{N}_{1}
\right) = \frac{1}{2} \Bigg(
\tr{\bs{\Sigma}_{1}^{-1} \bs{\Sigma}_{0}}
+ (\bs{\mu}_{1} - \bs{\mu}_{0})^{T} \bs{\Sigma}_{1}^{-1} 
\\
(\bs{\mu}_{1} - \bs{\mu}_{0}) - n 
+ \log{\frac{|\bs{\Sigma}_{1}|}{|\bs{\Sigma}_{0}|}} \Bigg). 
\end{multline}
\end{lemma}

\section{Main results}\label{section:Main_results}
In this section, the main results, which are divided into three parts, will be presented. First, two general results (\cref{theorem:An_inherently_unstable_system_cant_be_stabilized_covertly,theorem:An_inherently_stable_system_can_be_controlled_covertly}). Then two converse results under minimal-information controller (\cref{theorem:An_inherently_stable_system_cant_be_stabilized_covertly_one_bit,theorem:Observing_the_one_bit_controller_output_through_AWGN_channel}) and two results under maximal-information controller (\cref{theorem:Achievable_gain_using_threshold_controller,theorem:Converse_conditional_distribution}). 
The complete proofs are given in \cref{section:Proofs}. 

\cref{theorem:An_inherently_unstable_system_cant_be_stabilized_covertly} states that an inherently unstable linear system, i.e., an AR(1)  system with $|a| > 1$, can not be covertly stabilized. On the other hand, \cref{theorem:An_inherently_stable_system_can_be_controlled_covertly} states that an inherently stable linear system, i.e., an AR(1) system with $|a| < 1$, can be covertly controlled. 
In \cref{theorem:An_inherently_unstable_system_cant_be_stabilized_covertly} 
the stabilization criteria in mind is the absolute $\gamma$-moment stability, i.e., $\expc{}{}{|X_{n}|^{\gamma}} \leq c < \infty$. Hence, the meaning of stabilizing an unstable system is to achieve a finite absolute $\gamma$-moment for the system state. 
Conversely, a system with $|a| < 1$ is already stable, thus, stabilization action is not needed. However, if Alice desires to alter or interfere with the system's operation, it is possible to do so while keeping Willie ignorant about those actions.
In \cref{theorem:An_inherently_stable_system_can_be_controlled_covertly} below, Alice desires to change the gain of the system, i.e., her goal is to change an AR(1) system with a gain of $a$, to an AR(1) system with a gain of $b$. 
\begin{theorem}\label{theorem:An_inherently_unstable_system_cant_be_stabilized_covertly}
Consider the linear stochastic system in \cref{eq:simple_linear_model}, with $|a| > 1$, an i.i.d.\ $Z_{n}$ and a control signal $U_{n}$ which keeps the system $\gamma$-moment stable. 
If Willie has a uniform bound on the system's $\gamma$-moment, i.e., $\expc{}{}{|X_{n}|^{\gamma}} \leq c < \infty$ and observes the system's output through a clean channel, he can achieve $1-\delta$-detection, i.e., identify the control operation.

If, in adition, $Z_{n}$ is a Gaussian noise, Willie needs to observe the system's output at the time sample which is at least
\[
n_{0}
\geq
\frac{\log{\frac{M\sqrt{a^{2} - 1}}{\sigma_{Z} Q^{-1}\left( \frac{1 - \frac{\delta}{2}}{2} \right)}}}
{\log{|a|}}.
\]


where $M \geq \sqrt[\gamma]{\frac{2c}{\delta}}$, to achieves $1-\delta$-detection.
\end{theorem}

\cref{theorem:An_inherently_unstable_system_cant_be_stabilized_covertly} asserts that any attempt to covertly stabilize an inherently unstable AR(1) system, under the conditions of the claim, will lead to a failure. Under Gaussian noise, the theorem also gives an easy handle on where should Willie observe the system for any specific $\delta$. In addition, the theorem can be easily generalized to a scenario in which Willie is observing the system's output through a noisy channel.
This will only force Willie to increase the length and position of his observation interval, but the heart of the result remains.
\begin{theorem}\label{theorem:An_inherently_stable_system_can_be_controlled_covertly}
Consider the linear stochastic system given in \cref{eq:simple_linear_model} for $|a| < 1$, in steady state, $Z_{n} \sim \mathcal{N}(0, \sigma_{Z}^{2})$ i.i.d.  and $U_{n} = (a-b) X_{n-1}$, where $0 < |a| < |b| < 1$ and $\sgn{a} = \sgn{b}$. 
If Willie observes the system's output through a clean channel, for a time window $n < \frac{2b}{b-a}$, knows $U_{n}$ and $b$ satisfies  $|b|<\sqrt{1 - (1-a^{2}) e^{-4 \epsilon^{2}}}$, 
then for any method of detection that Willie will use, Alice achieves an $\epsilon$-covertness for any $\epsilon > 0$.
\end{theorem}

In other words, \cref{theorem:An_inherently_stable_system_can_be_controlled_covertly} states that with no information constraint on Alice's behalf, an inherently stable AR(1) system can be covertly controlled, in the sense that Alice can change the system's gain, $a$, to a different one, $b$, without being detected by Willie.
Again, \cref{theorem:An_inherently_stable_system_can_be_controlled_covertly} can be easily generalized to the scenario in which Willie observes the system's output trough a noisy channel, which in turn enables Alice to use further away $b$'s then in the claim. However,  \cref{theorem:An_inherently_stable_system_can_be_controlled_covertly} gives us a look at the best case that Willie can have (a clean channel), and still assures us that Willie can not detect Alice's control actions. 


The proof of \cref{theorem:An_inherently_unstable_system_cant_be_stabilized_covertly} is based on the fact that an AR(1) process with $|a| > 1$ is unstable and its moments diverge \cite{Time_Seires_Analysis}, thus, by observing the system's output far enough, Willie can make a calculated decision if Alice is stabilizing her system or not. 
However, the proof of  \cref{theorem:An_inherently_stable_system_can_be_controlled_covertly} 
is more intricate. It is based on (from Alice's side) creating a system for which the distribution of an $n$-tuple is very close to the original one in terms of divergence. By doing so, Willie would not be able to distinguish between the two distributions. 
\subsection{Converse results under a minimal-information controller}\label{subsetion:Converse_results_under_a_minimal-information_controller}
Consider now a setting in which Willie observes the controller's output trough a noisy channel. 
The following theorem introduces the case in which Willie observes Alice's one-bit controller's output through a noisy channel. 
It asserts that Willie can achieve $1-\delta$-detection, for any $\delta > 0$, as the number of observations increases.

This scenario is realistic in cases Willie dominates the channel Alice might have between her and the system. In this cases, it is also probable that an observer will know his channel characteristics, but not always have the resources to observe the controller in each time sample in which the system operates. Yet, as \Cref{theorem:Observing_the_one_bit_controller_output_through_AWGN_channel} asserts, if Willie observes the controller's output for a time period large enough, Willie can detect Alice's control actions and keep the sum of his error probabilities arbitrary small. Moreover, \Cref{theorem:An_inherently_stable_system_cant_be_stabilized_covertly_one_bit} adds another dimension to these results, stating that even in the case in which Willie observes the system's output, an inherently stable AR(1) cannot be covertly controlled by Alice if she uses the one-bit controller. Note, however, that these two results are under a minimal information controller, limiting Alice's ability to gather enough information before each control decision. The results in the next subsection will show that with more information, Alice can do better and stay undetected.
\begin{theorem}\label{theorem:Observing_the_one_bit_controller_output_through_AWGN_channel}
Consider an AR(1) system with $a \neq 2$, which Alice stabilizes using the one-bit controller of \cref{definition:One_bit_controller}. Willie observes the control signal through an AWGN for $K$ consecutive samples. The noise variance in Willie's channel, $\sigma_{v}^{2}$, is known to Willie. 

For any $\delta >0$, Willie can achieve $1-\delta$-detection as long as the observation window satisfies $K > K_{0}(\delta)$ for 
\[
\begin{split}
K_{0}(\delta)
=&
\frac{4}{\delta \cdot \text{SNR}^{2}} \left( 1 + \sqrt{1 + 2 \cdot \text{SNR}}\right)^{2},
\end{split}
\]

where $\text{SNR} \triangleq \frac{E_{U}}{\sigma_{v}^{2}}$, and $E_{U} = \left( \frac{aB}{2-a} \right)^{2}$ is the average energy of the one-bit controller. 
\end{theorem}
Very roughly speaking, the proof is based on the fact that the one--bit controller is coarse, using strong control actions, and hence is detectable with energy detection. 
\begin{theorem}
\label{theorem:An_inherently_stable_system_cant_be_stabilized_covertly_one_bit}
Consider an AR(1) system with $|a| < 1$, and $U_{n}$ as the one-bit controller. 
If Willie observes the system's output through a clean channel for at least
\small
\[
K_{0}(\delta)
=
\frac{1}{E_{U}^{2}} \left( \sqrt{\frac{m_{Z}(4) - \sigma_{Z}^{4} + 4E_{U} \sigma_{Z}^{2}}{\delta/2}}
+
\sqrt{\frac{m_{Z}(4) - \sigma_{Z}^{4}}{\delta/2}}
\right)^{2}
\]
\normalsize

time samples, where $m_{Z}(4) \triangleq \expc{}{}{Z^{4}}$, $\sigma_{Z}^{2} \triangleq \var{}{}{Z} = \expc{}{}{Z^{2}}$,
$E_{U} \triangleq \left( \frac{aB}{2-a} \right)^{2}$ and $B$ is some deterministic number for which $|Z| \leq B$, 
he can decide with high confidence if Alice is controlling the system, i.e., for any $\delta > 0$, Willie achieves $1-\delta$-detection.
\end{theorem}


To make the expression of the minimal observation window, $K_{0}(\delta)$, more compact, Willie can take,
\[
K_{0}(\delta)
\leq
\tilde{K}_{0}(\delta)
=
\frac{4}{E_{U}^{2}} \frac{m_{Z}(4) - \sigma_{Z}^{4} + 4E_{U} \sigma_{Z}^{2}}{\delta/2},
\]

and still achieve $1-\delta$-detection for any $K > \tilde{K}_{0}(\delta)$. 
\subsection{Results under a maximal-information controller}
\label{subsetion:Results_under_a_maximal-information_controller}
Finally, we turn to a maximum--information controller, and Alice's ability to covertly control under it. \cref{theorem:Achievable_gain_using_threshold_controller} introduces an achievable range of gains of an AR(1) system, for which Alice can use the threshold controller while staying undetected by Willie.  I.e., the claim shows that Alice can achieve $\varepsilon$-covertness for any $\varepsilon > 0$. It is restricted, however, to the case of one control action in the system's operation time, and for the case of inherently stable AR(1) system, i.e., $|a| < 1$ and under steady-state conditions. On the other hand, Willie is observing the system's output through a clean channel, for the whole of the system's operation time, and he is unrestricted in terms of complexity or strategy used. 
\begin{theorem}\label{theorem:Achievable_gain_using_threshold_controller}
Consider a Gaussian AR(1) system with $|a| < 1$. Alice is using the threshold controller, and Willie is observing to the system's output through a clean channel. 
If the system is being controlled by resetting at one time sample $1 \leq \tau \leq N$, i.e., $X_{\tau + 1} = Z_{\tau + 1}$, and the system's gain, $a$, satisfies $|a| 
\leq 
\sqrt{1 - 2^{-4 \varepsilon^{2}}}$, then for any method of detection that Willie will use,  Alice achieves $\varepsilon$-covertness for any $\varepsilon > 0$. 
\end{theorem}
The proof is based on computing the output distributions under both hypotheses. Note, however, that the distribution under the control hypothesis depends on $\tau$, whose distribution is unknown. The result uses distribution--independent bounds on the relevant divergence.

Next, we introduce a converse result to \cref{theorem:Achievable_gain_using_threshold_controller}. 
\cref{theorem:Converse_conditional_distribution} asserts that Willie can achieve an $1-\delta$-detection for any $\delta > 0$, as long as the gain of the system is considerably close to one. Note, however, that it gives Willie extra power, knowing the time Alice might activate the control. On the other hand, Willie observes the system only at that time.
\begin{theorem}\label{theorem:Converse_conditional_distribution}
Consider a Gaussian AR(1) system. Alice is using the threshold controller, and Willie is observing the system's output through a clean channel. 
If Willie knows the one time sample in which the system is being resets $1 \leq \tau \leq N$, i.e., $X_{\tau + 1} = Z_{\tau + 1}$, and the system's gain, $a$, satisfies
\[
|a|
\geq
\sqrt{1 - \frac{\left( Q^{-1}\left( \frac{1 - \delta /2}{2} \right) \right)^{2}}{2 \log{\frac{2}{\delta}}}},
\]

then, there exists a detection method in which Willie achieves an $1-\delta$-detection for any $\delta > 0$.
\end{theorem}
The proof is based on constructively suggesting a detection method. 

\newpage
\addcontentsline{toc}{section}{\numberline{}References}
\bibliographystyle{IEEEtran}
\bibliography{bibliography}

\newpage
\section{Proofs}\label{section:Proofs}
In this section, the proofs of the theorems given in \cref{section:Main_results} and their supporting claims  will be presented. 

\subsection{Proofs of \cref{lemma:PDF_of_n_samples_AR_1_process,lemma:Divergence_of_gaussian_vector}}
For completeness we give here the proofs of \cref{lemma:PDF_of_n_samples_AR_1_process,lemma:Divergence_of_gaussian_vector}.
\begin{proof}(\cref{lemma:PDF_of_n_samples_AR_1_process})
One can see that 
\cref{claim:System_state} can be represented as a vectors multiplication, with $U_{n} = 0$ . Thus, by concatenation we have, 
\begin{equation*}
\begin{split}
\textbf{X}^{(n)}
=&
\left(
\begin{matrix}
X_{1} \\
X_{2} \\
X_{3} \\
\vdots \\
X_{n}
\end{matrix}
\right)
\\=&
\underbrace{
\sigma_{Z}
\left(
\begin{matrix}
1 & 0 & 0 & \cdots & 0 \\
a & 1 & 0 & \cdots & 0 \\
a^{2} & a & 1 & \cdots & 0 \\
\vdots & \vdots & \vdots & \ddots & \vdots \\
a^{n-1} & a^{n-2} & a^{n-3} & \cdots & 1 \\
\end{matrix}
\right)}_{\textbf{A}}
\underbrace{
\frac{1}{\sigma_{Z}}
\left(
\begin{matrix}
Z_{1} \\
Z_{2} \\
Z_{3} \\
\vdots \\
Z_{n}
\end{matrix}
\right)}_{\tilde{\textbf{Z}}}
=
\textbf{A} \tilde{\textbf{Z}},
\end{split}
\end{equation*}

where $\tilde{\textbf{Z}} \sim \mathcal{N}(\textbf{0}, \textbf{I}_{n})$, and $\textbf{I}_{n}$ is the identity matrix of size $n$. 
Hence, we have 
$\textbf{X}^{(n)} \sim \mathcal{N}(\textbf{0}, \bs{\Sigma})$, 
when, 
\begin{equation*}
\begin{split}
\mathbb{E} \left[ \textbf{X}^{(n)} \right]
=&
\textbf{A} \mathbb{E} \left[ \tilde{\textbf{Z}} \right]
=
\textbf{0},
\\
Cov \left[ \textbf{X}^{(n)} \right]
=&
\mathbb{E} \left[ \textbf{X}^{(n)} \textbf{X}^{(n) T} \right]
=
\textbf{A} \mathbb{E} \left[ \tilde{\textbf{Z}} \tilde{\textbf{Z}}^{T} \right] \textbf{A}^{T}
=
\textbf{A} \textbf{A}^{T}.
\end{split}
\end{equation*}

By properties of the rank of a matrix, we have, 
$\operatorname*{rank}(\bs{\Sigma})
=
\operatorname*{rank}(\textbf{A} \textbf{A}^{T}) 
= 
\operatorname*{rank}(\textbf{A}) 
= 
n$. 
Therefore, $\bs{\Sigma}$ 
is invertible, and \cref{eq:PDF_of_n_samples_AR_1_process} 
is well-defined. In addition, $[\textbf{A}]_{i,j} = \sigma_{Z} a^{i-j} u(i-j)$, where $u(\cdot)$ is the discrete step function, hence we have,
\begin{equation*}
\begin{split}
[\bs{\Sigma}]_{ij}
=&
[\textbf{A} \textbf{A}^{T}]_{ij}
\\=&
\sum_{k=1}^{n}{[\textbf{A}]_{ik} [\textbf{A}]_{jk}}
\\=&
\sigma_{Z}^{2} \sum_{k=1}^{n}{a^{i-k} u(i-k) a^{j-k} u(j-k)}
\\\stackrel{(a)}{=}&
\sigma_{Z}^{2} a^{i+j} \sum_{k=1}^{\min(i,j)}{a^{-2k}}
\\=&
\sigma_{Z}^{2} a^{i+j} \frac{a^{-2 \min(i,j)} - 1}{1 - a^{2}}
\\\stackrel{(b)}{=}&
\frac{\sigma_{Z}^{2}}{1 - a^{2}} \left( a^{|i-j|} - a^{i+j} \right), \quad \forall 1 \leq i,j \leq n,
\end{split}
\end{equation*}

where (a) is since $i \geq k$ and $j \geq k$, hence, 
$k \leq \min(i,j)$. (b) is since $i+j-2 \min(i,j) = |i-j|$.
\end{proof}

\begin{proof}(\cref{lemma:Divergence_of_gaussian_vector})
It is easy to see that, 
\begin{align*}
\log \left(
\frac{\mathcal{N}_{0}}{\mathcal{N}_{1}}
\right)
=&
\frac{1}{2}
\log \frac{|\bs{\Sigma}_{1}|}{|\bs{\Sigma}_{0}|} 
+
\frac{1}{2}(\textbf{x} - \bs{\mu}_{1})^{T} \bs{\Sigma}_{1}^{-1} (\textbf{x} - \bs{\mu}_{1}) 
\\& - 
\frac{1}{2} (\textbf{x} - \bs{\mu}_{0})^{T} \bs{\Sigma}_{0}^{-1} (\textbf{x} - \bs{\mu}_{0}). 
\end{align*}

Applying expectation with respect to $\mathcal{N}_{0}$ yields, 
\begin{align*}
\expc{\mathcal{N}_{0}}{}{\log \left(
\frac{\mathcal{N}_{0}}{\mathcal{N}_{1}}
\right)}
=&
\frac{1}{2} \expc{\mathcal{N}_{0}}{}{(\textbf{x} - \bs{\mu}_{1})^{T} \bs{\Sigma}_{1}^{-1} (\textbf{x} - \bs{\mu}_{1})}
\\&
- \frac{1}{2} \expc{\mathcal{N}_{0}}{}{(\textbf{x} - \bs{\mu}_{0})^{T} \bs{\Sigma}_{0}^{-1} (\textbf{x} - \bs{\mu}_{0})}
\\& +
\frac{1}{2}
\log \frac{|\bs{\Sigma}_{1}|}{|\bs{\Sigma}_{0}|}
\end{align*}

Solving each of the expectations above, 
\begin{align*}
\mathbb{E}_{\mathcal{N}_{0}} & \left[ 
(\textbf{x} - \bs{\mu}_{0})^{T} \bs{\Sigma}_{0}^{-1} (\textbf{x} - \bs{\mu}_{0})
\right]
\\=&
\expc{\mathcal{N}_{0}}{}{\tr{(\textbf{x} - \bs{\mu}_{0})^{T} \bs{\Sigma}_{0}^{-1} (\textbf{x} - \bs{\mu}_{0})}}
\\=&
\expc{\mathcal{N}_{0}}{}{\tr{\bs{\Sigma}_{0}^{-1} (\textbf{x} - \bs{\mu}_{0})(\textbf{x} - \bs{\mu}_{0})^{T}}}
\\=&
\tr{\bs{\Sigma}_{0}^{-1} \expc{\mathcal{N}_{0}}{}{ (\textbf{x} - \bs{\mu}_{0})(\textbf{x} - \bs{\mu}_{0})^{T} }}
\\=&
\tr{\bs{\Sigma}_{0}^{-1} \bs{\Sigma}_{0}}
\\=&
n, 
\end{align*}

second, 
\begin{align*}
\mathbb{E}_{\mathcal{N}_{0}} & \left[ 
(\textbf{x} - \bs{\mu}_{1})^{T} \bs{\Sigma}_{1}^{-1} (\textbf{x} - \bs{\mu}_{1})
\right]
\\=&
\tr{\bs{\Sigma}_{1}^{-1} \expc{\mathcal{N}_{0}}{}{ (\textbf{x} - \bs{\mu}_{1})(\textbf{x} - \bs{\mu}_{1})^{T}}}
\\\stackrel{(a)}{=}&
\tr{\bs{\Sigma}_{1}^{-1} (\bs{\Sigma_{0}} + (\bs{\mu}_{1} - \bs{\mu}_{0}) (\bs{\mu}_{1} - \bs{\mu}_{0})^{T})}
\\=&
\tr{\bs{\Sigma}_{1}^{-1} \bs{\Sigma_{0}}} + (\bs{\mu}_{1} - \bs{\mu}_{0})^{T} \bs{\Sigma}_{1}^{-1} (\bs{\mu}_{1} - \bs{\mu}_{0}), 
\end{align*}

where (a) is since, 
\begin{align*}
\mathbb{E}_{\mathcal{N}_{0}} & \left[ 
(\textbf{x} - \bs{\mu}_{1})(\textbf{x} - \bs{\mu}_{1})^{T}
\right]
\\=&
\expc{\mathcal{N}_{0}}{}{(\textbf{x} - \bs{\mu}_{0} +\bs{\mu}_{0} - \bs{\mu}_{1})  (\textbf{x} - \bs{\mu}_{0} +\bs{\mu}_{0} - \bs{\mu}_{1})^{T}}
\\=&
\expc{\mathcal{N}_{0}}{}{(\textbf{x} - \bs{\mu}_{0}) (\textbf{x} - \bs{\mu}_{0})^{T}}
\\&+
\expc{\mathcal{N}_{0}}{}{(\textbf{x} - \bs{\mu}_{0})} (\bs{\mu}_{0} - \bs{\mu}_{1})^{T} 
\\&+
(\bs{\mu}_{0} - \bs{\mu}_{1}) \expc{\mathcal{N}_{0}}{}{(\textbf{x} - \bs{\mu}_{0})^{T}} 
\\&+
(\bs{\mu}_{0} - \bs{\mu}_{1}) (\bs{\mu}_{0} - \bs{\mu}_{1})^{T}
\\=&
\bs{\Sigma_{0}} + (\bs{\mu}_{1} - \bs{\mu}_{0}) (\bs{\mu}_{1} - \bs{\mu}_{0})^{T}.
\end{align*}
\end{proof}

\subsection{Proof of \cref{theorem:An_inherently_unstable_system_cant_be_stabilized_covertly}}

\begin{proof}(\cref{theorem:An_inherently_unstable_system_cant_be_stabilized_covertly})
An AR(1) process with a gain $|a| > 1$ is not stable, in the sense that the second moment of the process goes to infinity with time, i.e., $\expc{}{}{X_{n}^{2}} \xrightarrow[n \to \infty]{} \infty$ \cite{Time_Seires_Analysis}. 
Therefore, if Willie observes the system's output through a clean channel, at a far enough time sample, he can decide with high confidence if Alice is stabilizing the system or not, namely, if the system state stays in some bounded region, or tents to infinity, respectively. 


 
The detection method Willie uses is as follows: Willie observes the system's output at time $n_{0}$ and check if its absolute value is bigger or smaller than some deterministic number $M$, i.e., Willie observes $X_{n_{0}}$ and compares $|X_{n_{0}}|$ with $M$. 

This analysis can be formulated as the following hypothesis testing problem,  
\begin{align*}
\begin{matrix}
\mathcal{H}_{0}: &
X_{n} = a X_{n-1} + Z_{n}, & 
n \geq 1, 
\\
\mathcal{H}_{1}: &
X_{n} = a X_{n-1} + Z_{n} - U_{n}, & 
n \geq 1,
\end{matrix}
\end{align*}

where $U_{n}$ is such that $\expc{}{}{|X_{n_{0}}|^{\gamma}} \leq c$. 
Hence, we bound the miss detection probability as follows, 
\begin{align*}
\beta 
=&
\prob{|X_{n_{0}}| > M | \mathcal{H}_{1}}
\\ \stackrel{(a)}{\leq}&
\frac{\expc{}{}{|X_{n_{0}}|^{\gamma} | \mathcal{H}_{1}}}{M^{\gamma}}
\\ \stackrel{(c)}{\leq}&
\frac{c}{M^{\gamma}}
\\ \stackrel{(d)}{\leq}&
\frac{\delta}{2},
\end{align*}

where (a) is due to Markov's inequality, (b) since $U_{n}$ stabilizes the system in the sense that $\expc{}{}{|X_{n_{0}}|^{\gamma} | \mathcal{H}_{1}} \leq c < \infty$. 
(d) is when Willie sets $M \geq \sqrt[\gamma]{\frac{2c}{\delta}}$ to get $\beta \leq \frac{\delta}{2}$. 

On the other hand, if $\mathcal{H}_{0}$ is true, than $|X_{n_{0}}|$ should be big with high probability. 
Thus, the probability of Willie to falsely decide that Alice is stabilizing the system when she is not, is when $|X_{n_{0}}|$ is smaller than some big enough constant $M$, hence, we evaluate the false alarm probability as follows, 
\begin{align*}
\alpha
=&
\prob{|X_{n_{0}}| \leq M | \mathcal{H}_{0}}
\\=&
\prob{\left\vert \sum_{k=1}^{n_{0}}{a^{n_{0}-k} Z_{k}} \right\vert < M}
\\\stackrel{(a)}{=}&
\prob{- \frac{M}{|a|^{n_{0}}} < \sum_{k=1}^{n_{0}}{a^{-k} Z_{k}} < \frac{M}{|a|^{n_{0}}}}.
\end{align*}

For $|a| > 1$ the sum $\tilde{Z} = \sum_{k=1}^{n_{0}}{a^{-k} Z_{k}}$ has a variance of $\var{}{}{\tilde{Z}} = \sigma_{Z}^{2} \frac{1 - a^{-2n_{0}}}{a^{2} - 1}$, which is finite and nonzero for any $n_{0}$, even for $n_{0} \to \infty$. 
On the other hand, $\frac{M}{|a|^{n_{0}}} \xrightarrow[n_{0} \to \infty]{} 0$, hence $\alpha \xrightarrow[n_{0} \to \infty]{} 0$. 

In the special case in which $Z_{k} \sim \mathcal{N}(0,\sigma_{Z}^{2})$, we have 
\begin{align*}
\alpha
=&
\prob{- \frac{M}{|a|^{n_{0}}} < \sum_{k=1}^{n_{0}}{a^{-k} Z_{k}} < \frac{M}{|a|^{n_{0}}}}
\\=&
1 - 2Q\left( \frac{M\sqrt{a^{2} - 1}}{\sigma_{Z} \sqrt{a^{2n} - 1}} \right)
\\ \stackrel{(a)}{\leq}&
\frac{\delta}{2}
\end{align*}

where (a) is by applying the detection constraint. Rearranging terms yields the requirement 
\begin{align*}
n_{0}
\geq&
\log_{|a|}{\frac{M\sqrt{a^{2} - 1}}{\sigma_{Z} Q^{-1}\left( \frac{1 - \frac{\delta}{2}}{2} \right)}}
\\=&
\frac{\log{\frac{M\sqrt{a^{2} - 1}}{\sigma_{Z} Q^{-1}\left( \frac{1 - \frac{\delta}{2}}{2} \right)}}}
{\log{|a|}}.
\end{align*}

\end{proof}

\subsection{Proof of \cref{theorem:An_inherently_stable_system_can_be_controlled_covertly}}
\label{subsubsection:An_inherently_stable_system_can_be_stabilized_covertly}
First, we prove the next supporting claim.
\begin{claim}\label{claim:trace_B_inv_A}
Consider two $n$-tuple Gaussian AR(1) processes in steady state. 
The first with a gain of $a$ and a covariance matrix $\bs{\Sigma}_{0}$, and the other with a gain of $b$ and a covariance matrix $\bs{\Sigma}_{1}$. 
For $|a|, |b| < 1$, we have
\begin{align*}
\tr{\bs{\Sigma}_{1}^{-1} \bs{\Sigma}_{0}}
=&
\frac{(n-2) b^{2} - 2(n-1) ab + n}{1 - a^{2}}.
\end{align*}
\end{claim}

\begin{proof}
For $|a|, |b| < 1$, and in steady state, by  \cref{lemma:PDF_of_n_samples_AR_1_process} and  \cref{corol:Stable_AR(1)_Covariance_in_SS}, 
$[\bs{\Sigma_{0}}]_{i,j} = \frac{\sigma_{Z}^{2}}{1 - a^{2}} a^{|i-j|}$ and 
$[\bs{\Sigma_{1}}]_{i,j} = \frac{\sigma_{Z}^{2}}{1 - b^{2}} b^{|i-j|}$ for $1 \leq i,j \leq n$. 
Denote $[\textbf{A}]_{i,j} = a^{|i-j|}$ and 
$[\textbf{B}]_{i,j} = b^{|i-j|}$ for $1 \leq i,j \leq n$. 
To show that 
\begin{align*}
[\textbf{A}^{-1}]_{i,j}
=&
\frac{
\left( 1 + a^{2} 1_{2 \leq i \leq n-1} \right) 1_{i=j} - a \left( 1_{i = j-1} + 1_{i = j+1} \right)
}{1 - a^{2}},
\end{align*}

we check by definition that $\textbf{A}^{-1} \textbf{A} = \textbf{I}_{n}$. 
First, 
\begin{align*}
[\textbf{A}^{-1} \textbf{A}]_{i,j}
=&
\sum_{k=1}^{n}{[\textbf{A}]_{i,k} [\textbf{A}^{-1}]_{k,j}}
\\=&
\frac{1}{1 - a^{2}}
\sum_{k=1}^{n}{\left[
\left( 1 + a^{2} 1_{2 \leq k \leq n-1} \right) 1_{k=j} a^{|i-k|} 
\right.}
\\& \left.
- a \left( 1_{k = j-1} + 1_{k = j+1} \right) a^{|i-k|} \right]
\\=&
\frac{1}{1 - a^{2}}
\left[
\left( 1 + a^{2} 1_{2 \leq j \leq n-1} \right) a^{|i-j|} \right.
\\& \left.
- a \left( a^{|i-j + 1|} 1_{j \neq 1} + a^{|i-j - 1|} 1_{j \neq n} \right) \right],
\end{align*}

for $i \geq j + 1$, $2 \leq i \leq n$ and $1 \leq j \leq n-1$, thus
\begin{align*}
[\textbf{A}^{-1} \textbf{A}]_{i,j}
=&
\frac{1}{1 - a^{2}}
\left[
\left( 1 + a^{2} 1_{2 \leq j \leq n-1} \right) a^{i-j} \right.
\\& \left.
- a \left( a^{i-j + 1} 1_{j \neq 1} + a^{i-j - 1} 1_{j \neq n} \right) \right]
\\=&
\frac{1}{1 - a^{2}}
\left[
a^{i-j} 1_{j = n} 
+
a^{i-j+2} 1_{j = n} 
\right]
= 0,
\end{align*}

Similarly for $j \geq i +1$, we have $[\textbf{A}^{-1} \textbf{A}]_{i,j} = 0$. For $i = j$, 
\begin{align*}
[\textbf{A}^{-1} \textbf{A}]_{i,i}
=&
\frac{1}{1 - a^{2}}
\left[
\left( 1 + a^{2} 1_{2 \leq j \leq n-1} \right) \right.
\\& \left.
- a \left( a 1_{j \neq 1} + a 1_{j \neq n} \right) \right]
\\=&
\frac{1}{1 - a^{2}}
\left[
1 - a^{2}
\right]
= 1.
\end{align*}

Now, by changing $a$ to $b$ in $[\textbf{A}^{-1}]_{i,j}$,  $[\textbf{B}^{-1}]_{i,j}$ is obtained. 
Next, we evaluate the following trace 
\begin{align*}
\tr{\textbf{B}^{-1} \textbf{A}}
=&
\sum_{m=1}^{n}{[\textbf{B}^{-1} \textbf{A}]_{m,m}}
\\=&
\sum_{m=1}^{n}{
\sum_{k=1}^{n}{
[\textbf{B}^{-1}]_{m,k} [\textbf{A}]_{k,m}}}
\\=&
\frac{1}{1 - b^{2}}
\sum_{m=1}^{n}{
\sum_{k=1}^{n}{
\left[ 
-b a^{|k-m|} (1_{k=m-1} + 1_{k=m+1})
\right.}}
\\& \left.
+ b^{2} a^{|k-m|} 1_{2 \leq m \leq n-1} 1_{k=m}
+ a^{|k-m|} 1_{k=m}
\right]
\\=&
\frac{1}{1 - b^{2}}
\sum_{m=1}^{n}{
\left[ 
-ab (1_{m \neq 1} + 1_{m \neq n})
\right.}
\\&+ \left.
b^{2} 1_{2 \leq m \leq n-1}
+ 1
\right]
\\=&
\frac{(n-2)b^{2} - 2(n-1) ab + n}{1 - b^{2}}.
\end{align*}

Thus, 
\begin{align*}
\tr{\bs{\Sigma}_{1}^{-1} \bs{\Sigma}_{0}}
=&
\frac{1 - b^{2}}{1 - a^{2}}
\tr{\textbf{B}^{-1} \textbf{A}}
\\=&
\frac{(n-2) b^{2} - 2(n-1) ab + n}{1 - a^{2}}.
\end{align*}
\end{proof}

\begin{proof}(\cref{theorem:An_inherently_stable_system_can_be_controlled_covertly})
Alice's goal is to convert a given AR(1) system with a gain of $a$, to another AR(1) system with a gain of $b$. 
In Willie's side, this case can be formulated as the following hypothesis testing problem,  
\begin{align*}
\begin{matrix}
\mathcal{H}_{0}: &
X_{n} = a X_{n-1} + Z_{n}, & 
n \geq 1,
\\
\mathcal{H}_{1}: &
X_{n} = b X_{n-1} + Z_{n}, & 
n \geq 1 .
\end{matrix}
\end{align*}

Since $Z_{n} \sim \mathcal{N}(0, \sigma_{Z}^{2})$ i.i.d., then $\textbf{X}^{(n)}| \mathcal{H}_{0} \sim \mathcal{N} \left( \textbf{0}, \bs{\Sigma_{0}} \right)$ and
$\textbf{X}^{(n)}| \mathcal{H}_{1} \sim \mathcal{N} \left( \textbf{0}, \bs{\Sigma_{1}} \right)$, where $[\bs{\Sigma_{0}}]_{i,j} = \frac{\sigma_{Z}^{2}}{1 - a^{2}} a^{|i-j|}$ and $[\bs{\Sigma_{1}}]_{i,j} = \frac{\sigma_{Z}^{2}}{1 - b^{2}} b^{|i-j|}$ for $1 \leq i,j \leq n$, respectively 
(see \cref{lemma:PDF_of_n_samples_AR_1_process} and \cref{corol:Stable_AR(1)_Covariance_in_SS}). 

Hence, 
\[
\begin{split}
\alpha + \beta
\stackrel{(a)}{=}&
1 - 
\mathcal{V}_{T}\left( f_{\textbf{X}| \mathcal{H}_{0}}, f_{\textbf{X}| \mathcal{H}_{1}} \right)
\\ \stackrel{(b)}{\geq}&
1 - 
\sqrt{\frac{1}{2} \dive{KL}{f_{\textbf{X}| \mathcal{H}_{0}}}{f_{\textbf{X}| \mathcal{H}_{1}}}}
\\ \stackrel{(c)}{=}&
1 - 
\sqrt{\frac{1}{4} \left( \tr{\bs{\Sigma}_{1}^{-1} \bs{\Sigma}_{0}} - n + \log{\left( \frac{\left\vert  \bs{\Sigma}_{1} \right\vert}{\left\vert  \bs{\Sigma}_{0} \right\vert} \right)} \right)}
\\ \stackrel{(d)}{>}&
1 - \frac{1}{2} 
\sqrt{\log{\left( \frac{1 - a^{2}}{1 - b^{2}} \right)}}
\\ \stackrel{(e)}{>}&
1 - 
\varepsilon,
\end{split}
\]

(a) 
is due to 
\cref{lemma:Connection_between_error_probabilities_to_total_variation}. (b) is due to 
\cref{lemma:Connection_between_divergence_to_total_variation}. 
(c) 
is by \cref{lemma:Divergence_of_gaussian_vector}. 
(d) is due to the assumption that $|b| > |a|$, $\sgn{a} = \sgn{b}$ and $n < \frac{2b}{b-a}$, hence by \cref{claim:trace_B_inv_A} 
\begin{align*}
\tr{\bs{\Sigma}_{1}^{-1} \bs{\Sigma}_{0}} - n
=&
\frac{(b-a)^{2}}{1 - a^{2}} n - 2b \frac{b-a}{1 - a^{2}}
\\=&
\frac{(b-a)^{2}}{1 - a^{2}}\left(
n  - \frac{2b}{b-a}
\right)
\\<&
0.
\end{align*}

The logarithm is by the substitution of $\left\vert  \bs{\Sigma}_{0} \right\vert = \frac{\sigma_{Z}^{2n}}{1 - a^{2}}$ and $\left\vert  \bs{\Sigma}_{1} \right\vert = \frac{\sigma_{Z}^{2n}}{1 - b^{2}}$ (see the proof of \cref{claim:Divergence_of_controlled_and_uncontrolled_AR_1_of_stable_system_in_steady_state}). 
(e) is by applying the covertness criterion. 
This results in 
\[
\varepsilon
>
\frac{1}{2} \sqrt{\log{\left( \frac{1 - a^{2}}{1 - b^{2}} \right)}},
\]

or, alternatively, the requirement 
\[
|b| 
<
\sqrt{1 - (1-a^{2}) e^{-4 \epsilon^{2}}}.
\]
\end{proof}

\subsection{Proof of \cref{theorem:Observing_the_one_bit_controller_output_through_AWGN_channel}}




In this subsection, the proof of  \cref{theorem:Observing_the_one_bit_controller_output_through_AWGN_channel} 
and its supporting claims will be presented. We start with three supporting claims.
\begin{claim}\label{claim:System_state}
The linear stochastic system shown in \cref{eq:simple_linear_model}, initializes at $n = 0$ (thus $X_{0} = 0$), 
can be represented for an arbitrary controller $U_{n}$ as follows,
\begin{align}\label{eq:System_state}
\begin{split}
X_{k} = \sum_{m = 1}^{k}{a^{k - m} \left( Z_{m} -  U_{m} \right)}, \quad  \forall k \in \mathbb{N}, \ a \in \mathbb{R},
\end{split}
\end{align}

where $X_{k}$ is the system state at time $k$.
\end{claim}

\begin{proof}
The proof is by induction. Clearly, for $k = 1$ the claim holds since,
\begin{align*}
\begin{split}
X_{1}
=&
\sum_{m = 1}^{1}{a^{k - m} \left( Z_{m} -  U_{m} \right)}
\\=&
a^{0} (Z_{1} -  U_{1})
\\=&
Z_{1}, 
\end{split}
\end{align*}

where $U_{1} = 0$. 
Now, for $k+1$, we have
\begin{align*}
\begin{split}
\sum_{m = 1}^{k+1}{a^{k + 1 - m} \left( Z_{m} -  U_{m} \right)}
=&
\sum_{m = 1}^{k}{a^{k + 1 - m} \left( Z_{m} -  U_{m} \right)} 
\\&
+ Z_{k + 1} - U_{k + 1}
\\=&
a \left[\sum_{m = 1}^{k}{a^{k - m} \left( Z_{m} -  U_{m} \right)}\right] 
\\&
+ Z_{k + 1} - U_{k + 1}
\\=&
a X_{k} + Z_{k + 1} - U_{k + 1}
\triangleq
X_{k + 1} .
\end{split}
\end{align*}

Hence, by mathematical induction \cref{eq:System_state} holds. 
\end{proof}

\begin{claim}\label{claim:One_bit_controller_representation}
The \emph{One-bit controller}, shown in \cref{eq:One_bit_controller}, can be written as,
\begin{align}\label{eq:One_bit_controller_representation}
\begin{split}
U_{n}
=
\left[
\frac{a}{2} \frac{B}{1 - a/2} + \left(\frac{a}{2}\right)^{n-1} \left( C_{1} - \frac{B}{1 - a/2}\right)
\right] \operatorname*{sgn}(X_{n-1}),
\end{split}
\end{align}

where 
$C_{1} \geq \frac{B}{1 - a/2}$ 
is the first element in the series $C_{n}$, 
$B$ is the bound of the noise $Z_{n}$, .i.e., $|Z_{n}| \leq B, \ \forall n \geq 1$ and  
$a \neq 2$ is the gain of the system shown in  \cref{eq:simple_linear_model}.
\end{claim}

\begin{proof}
By \cref{eq:One_bit_controller}, one has to show that the following holds,
\begin{align}\label{eq:One_bit_controller_C_n}
\begin{split}
C_{n-1}
=
\frac{B}{1 - a/2} + \left(\frac{a}{2}\right)^{n-2} \left( C_{1} - \frac{B}{1 - a/2}\right).
\end{split}
\end{align}

First, we will prove the following by induction,
\begin{align}\label{eq:One_bit_controller_C_n-k}
\begin{split}
C_{n}
=
\frac{B}{1 - a/2}
+ \left(\frac{a}{2}\right)^{k} \left( C_{n-k} - \frac{B}{1 - a/2} \right),
\end{split}
\end{align}

where $k$ is some shift of the series $C_{n}$. For instance, if $k=1$ we have $C_{n} = (a/2) C_{n-1} + B$, which is by definition. Assume \cref{eq:One_bit_controller_C_n-k} holds for 
$k=r$, then for $k= r+1$,
\begin{align*}
\begin{split}
C_{n}
=&
\frac{B}{1 - a/2}
+ \left(\frac{a}{2}\right)^{r} \left( C_{n-r} - \frac{B}{1 - a/2} \right)
\\\stackrel{(a)}{=}&
\frac{B}{1 - a/2}
+ \left(\frac{a}{2}\right)^{r} \left( \left( \frac{a}{2} C_{n-r-1} + B \right) - \frac{B}{1 - a/2} \right)
\\=&
\frac{B}{1 - a/2}
+ \left(\frac{a}{2}\right)^{r} \left( \frac{a}{2} C_{n-r-1} - \frac{a}{2} \frac{B}{1 - a/2} \right)
\\=&
\frac{B}{1 - a/2}
+ \left(\frac{a}{2}\right)^{r+1} \left( C_{n-(r+1)} - \frac{B}{1 - a/2} \right),
\end{split}
\end{align*}

where (a) is by the relation: $C_{n} = (a/2) C_{n-1} + B$. 
Therefore, \cref{eq:One_bit_controller_C_n-k} holds $\forall 1 \leq k < n$. 
Let us substitute to \cref{eq:One_bit_controller_C_n-k}, $k = n-1$,
\begin{align}\label{eq:One_bit_controller_C_n_proved}
\begin{split}
C_{n}
=&
\frac{B}{1 - a/2}
+ \left(\frac{a}{2}\right)^{n-1} \left( C_{n-(n-1)} - \frac{B}{1 - a/2} \right)
\\=&
\frac{B}{1 - a/2}
+ \left(\frac{a}{2}\right)^{n-1} \left( C_{1} - \frac{B}{1 - a/2} \right).
\end{split}
\end{align}
\end{proof}

\begin{coroc}
By setting $C_{1} = \frac{B}{1 - a/2}$,
\begin{align}\label{eq:One_bit_controller_representation_spcial_case}
\begin{split}
U_{n}
=
\frac{a}{2} \frac{B}{1 - a/2} \operatorname*{sgn}(X_{n-1}), \quad \forall n > 1,
\end{split}
\end{align}

since $C_{n}$ is monotonically decreasing to $\frac{B}{1 - a/2}$, hence,
\begin{align}\label{eq:One_bit_controller_representation_spcial_case_abs_value}
\begin{split}
|U_{n}|
=
\frac{a}{2} \frac{B}{1 - a/2}, \quad \forall n > 1,
\end{split}
\end{align}

since $a/2 < 1$. 
As $n \to \infty$, \crefrange{eq:One_bit_controller_representation_spcial_case}{eq:One_bit_controller_representation_spcial_case_abs_value} also holds regardless of the choice of $C_{1}$ (however, it has to be greater or equal to $\frac{B}{1 - a/2}$). 
Thus, one can deduce that $|U_{n}|$ 
converges. I.e., $|U_{n}|$ can not get arbitrary large.
\end{coroc}

Next, we give an upper and lower bounds for the energy of the control signal shown in \cref{eq:One_bit_controller}. 
\begin{claim}\label{claim:One_bit_controller_energy}
The energy of the controller in  \cref{eq:One_bit_controller}, is bounded by,
\begin{align}\label{eq:One_bit_controller_energy}
\begin{split}
\left( \frac{a B}{2 - a} \right)^{2} \leq E_{U} \leq \left( \frac{a}{2} C_{1} \right)^{2},
\end{split}
\end{align}
 
where 
$C_{n}$ is a deterministic, monotonically decreasing series, converging to: $\frac{B}{1 - a/2}$. Moreover, 
$C_{1}$ can be arbitrary number which sustains $C_{1} \geq \frac{B}{1- a/2}$, and $C_{n} = (a/2) C_{n-1} + B$, where $B$ is the noise bound, and $a \in \mathbb{R},  a \neq 2$ (see \cref{definition:One_bit_controller}).
\end{claim}

\begin{proof}
we have,
\begin{align*}
\begin{split}
E_{U} 
\triangleq&
\frac{1}{N}\sum_{n = 1}^{N}{U_{n}^{2}}
\\=&
\frac{1}{N}\sum_{n = 1}^{N}{\frac{a^{2}}{4} C_{n-1}^{2} (\operatorname*{sgn}(X_{n-1}))^{2}}
\\=&
\frac{a^{2}}{4} \frac{1}{N}\sum_{n = 1}^{N}{C_{n-1}^{2}}
\\\stackrel{(a)}{\leq}&
\frac{a^{2}}{4} \frac{1}{N}\sum_{n = 1}^{N}{C_{1}^{2}}
\\=&
\left( \frac{a}{2} C_{1} \right)^{2},
\end{split}
\end{align*}

where (a) is since $C_{n}$ is a monotonically decreasing series. 

On the other hand,
\begin{align*}
\begin{split}
E_{U} 
=&
\frac{a^{2}}{4} \frac{1}{N}\sum_{n = 1}^{N}{C_{n-1}^{2}}
\\\stackrel{(a)}{\geq}&
\frac{a^{2}}{4} \frac{1}{N}\sum_{n = 1}^{N}{\left( \frac{B}{1 - a/2} \right)^{2}}
\\=&
\frac{a^{2}}{4} \left( \frac{B}{1 - a/2} \right)^{2}
\\=&
\left( \frac{a B}{2 - a} \right)^{2},
\end{split}
\end{align*}

where (a) is since $C_{n}$ is a monotonically decreasing series converging to: $\frac{B}{1 - a/2}$.
\end{proof}

\begin{coroc}\label{coroc:One_bit_controller_energy_steady_state}
Since $C_{1}$ can be arbitrary number which sustains $C_{1} \geq \frac{B}{1- a/2}$, 
for simplicity, we set $C_{1} = \frac{B}{1 - a/2}$ unless otherwise stated, thus,  
\begin{align}\label{eq:One_bit_controller_energy_steady_state}
\begin{split}
E_{U} = \left( \frac{a B}{2 - a} \right)^{2},
\end{split}
\end{align}

which yields a constant energy with respect to time.

Furthermore, from \cref{eq:One_bit_controller_representation_spcial_case_abs_value}, 
the above also holds as $N \to \infty$ regardless of the choice of $C_{1}$ (however, $C_{1}$ has to be greater or equal to $\frac{B}{1 - a/2}$). 
Hence, the average energy of \cref{eq:One_bit_controller} in steady state, is given by \cref{eq:One_bit_controller_energy_steady_state} regardless of the choice of $C_{1}$.
\end{coroc}

By claim 2 in \cite{bash_limits}, one can show that Alice cannot covertly control her system, as Willie's observations window increases. 
However, we will show a detailed proof for 
\cref{theorem:Observing_the_one_bit_controller_output_through_AWGN_channel}. 
  
\begin{proof}(\cref{theorem:Observing_the_one_bit_controller_output_through_AWGN_channel})
Willie observes Alice's control signal through an AWGN for $K$ consecutive samples. 
Willie's observation at time sample $n$ is simply $W_{n} = U_{n} + v_{n}$, where $U_{n}$ is Alice's control signal, $v_{n} \sim \mathcal{N}(0, \sigma_{v}^{2})$ is an i.i.d.\ noise and the noise variance, $\sigma_{v}^{2}$, is known to Willie. 
Since, $U_{n} \neq 0 \quad \forall n > 1$, i.e., Alice always controlling the system, then for any $K$ consecutive samples the average energy that Willie reads is: $E_{W} = \frac{1}{K}\sum_{n}^{}{(U_{n} + v_{n})^{2}} $. 
Hence, Willie will compare $E_{W}$ to the noise energy in his channel, i.e., $\sigma_{v}^{2}$. 


Willie has the following hypotheses
\begin{align*}
\begin{matrix}
\mathcal{H}_{0}: &
W_{n} = v_{n}, & 
n = k, \ldots, K + k - 1,
\\
\mathcal{H}_{1}: &
W_{n} = U_{n} + v_{n}, &
n = k, \ldots, K + k - 1.
\end{matrix}
\end{align*}

Thus, the mean and the variance of his energy, $E_{W}$, when $\mathcal{H}_{0}$ is true are,
\begin{align}\label{eq:Analytical_analysis:One_bit:Covert:Mean_H_0}
\begin{split}
\mathbb{E} \left[ E_{W}| \mathcal{H}_{0} \right]
=&
\sigma_{v}^{2},
\end{split}
\end{align}
\begin{align}\label{eq:Analytical_analysis:One_bit:Covert:Var_H_0}
\begin{split}
\var{}{}{ E_{W}| \mathcal{H}_{0}}
=&
\var{}{}{\frac{1}{K}\sum_{n}^{}{v_{n}^{2}}}
\\=&
\frac{1}{K^{2}}\sum_{n}^{}{ \left( \mathbb{E} \left[ v_{n}^{4}\right] - \mathbb{E}^{2} \left[ v_{n}^{2} \right] \right)} 
\\=&
\frac{2 \sigma_{v}^{4}}{K},
\end{split}
\end{align}

since $v_{n} \sim \mathcal{N}(0, \sigma_{v}^{2})$ and  i.i.d. Under the alternative hypothesis, the mean and the variance of $E_{W}$ are,
\begin{align}\label{eq:Analytical_analysis:One_bit:Covert:Mean_H_1}
\begin{split}
\mathbb{E} \left[ E_{W}| \mathcal{H}_{1} \right]
=&
\mathbb{E} \left[ \frac{1}{K}\sum_{n}^{}{(U_{n} + v_{n})^{2}} \right]
\\=&
\frac{1}{K} \sum_{n}^{}{\expc{}{}{U_{n}^{2}}}
+ 
\frac{2}{K}\sum_{n}^{}{\mathbb{E} \left[ U_{n} v_{n} \right]} 
\\&+
\frac{1}{K}\sum_{n}^{}{\mathbb{E} \left[ v_{n}^{2}\right]} 
\\=&
E_{U}
+
\sigma_{v}^{2},
\end{split}
\end{align}
\begin{align}\label{eq:Analytical_analysis:One_bit:Covert:Second_moment_H_1}
\begin{split}
\mathbb{E} \left[ E_{W}^{2}| \mathcal{H}_{1} \right]
=&
\mathbb{E} \left[ \frac{1}{K^{2}}\sum_{n,m}^{}{(U_{n} + v_{n})^{2} (U_{m} + v_{m})^{2}} \right]
\\=&
\frac{1}{K^{2}}\sum_{n,m}^{}{ \left( U_{n}^{2} U_{m}^{2} + 
2 U_{n}^{2} \mathbb{E} \left[ U_{m} v_{m}  \right] \right.}
\\&+ 
U_{n}^{2} \mathbb{E} \left[ v_{m}^{2}  \right] 
+
2 U_{m}^{2} \mathbb{E} \left[ U_{n} v_{n}  \right] 
\\&+ 
4 \mathbb{E} \left[ U_{n} v_{n} U_{m} v_{m} \right] 
+ 2 \mathbb{E} \left[ U_{n} v_{n} v_{m}^{2} \right] 
\\&+ 
U_{m}^{2} \mathbb{E} \left[ v_{n}^{2} \right] 
+
2 \mathbb{E} \left[U_{m} v_{m} v_{n}^{2} \right] 
\\&+  \left. 
\mathbb{E} \left[ v_{m}^{2} v_{n}^{2} \right]
\right)
\\=&
E_{U}^{2} + 2 E_{U} \sigma_{v}^{2} + 
\frac{4 E_{U} \sigma_{v}^{2} + (K+2) \sigma_{v}^{4}}{K},
\end{split}
\end{align}

since $U_{n}^{2}$ is deterministic. Using \crefrange{eq:Analytical_analysis:One_bit:Covert:Mean_H_1}{eq:Analytical_analysis:One_bit:Covert:Second_moment_H_1}, 
\begin{align}\label{eq:Analytical_analysis:One_bit:Covert:Var_H_1}
\begin{split}
\var{}{}{E_{W}| \mathcal{H}_{1}}
=&
\mathbb{E} \left[ E_{W}^{2} \right]
-
\mathbb{E}^{2} \left[ E_{W} \right]
\\=&
E_{U}^{2} + 2 E_{U} \sigma_{v}^{2} + \sigma_{v}^{4} + \frac{4 E_{U} \sigma_{v}^{2} + 2 \sigma_{v}^{4}}{K} 
\\&-
\left[ E_{U}^{2} + 2E_{U} \sigma_{v}^{2} + \sigma_{v}^{4} \right]
\\=&
\frac{4 E_{U} \sigma_{v}^{2} + 2 \sigma_{v}^{4}}{K},
\end{split}
\end{align}

where $E_{U} = \left( \frac{aB}{2-a}  \right)^{2}$  is the average energy of the one-bit controller (\cref{coroc:One_bit_controller_energy_steady_state}). 

Willie picks a threshold which we denote as $t$, and compares the value of $E_{W}$ to $\sigma_{v}^{2} + t$. Willie accepts $\mathcal{H}_{0}$ if  $E_{W} < \sigma_{v}^{2} + t$ and rejects it otherwise. Suppose that Willie desires a false alarm probability which is bounded by $\delta/2$, which is the probability
that $E_{W} \geq \sigma_{v}^{2} + t$ when $\mathcal{H}_{0}$ is true. By 
\crefrange{eq:Analytical_analysis:One_bit:Covert:Mean_H_0}{eq:Analytical_analysis:One_bit:Covert:Var_H_0}, using Chebyshev’s inequality, we have, 
\begin{align*}
\begin{split}
\alpha
=&
\mathbb{P}\{ E_{W} \geq \sigma_{v}^{2} + t| \mathcal{H}_{0} \}
\\\leq&
\mathbb{P}\{ |E_{W} - \sigma_{v}^{2}| \geq t| \mathcal{H}_{0} \}
\\\leq&
\frac{2 \sigma_{v}^{4}}{K t^{2}}.
\end{split}
\end{align*}

Thus, to obtain $\alpha \leq \delta/2$, Willie sets $t = \sqrt{\frac{4}{\delta}} \frac{\sigma_{v}^{2}}{\sqrt{K}}$. 
The probability of a miss detection, $\beta$, is the probability that $E_{W} < \sigma_{v}^{2} + t$ when $\mathcal{H}_{1}$ is true. 
By  \crefrange{eq:Analytical_analysis:One_bit:Covert:Mean_H_1}{eq:Analytical_analysis:One_bit:Covert:Var_H_1}, using Chebyshev’s inequality, we have,
\begin{align*}
\begin{split}
\beta
=&
\mathbb{P}\{ E_{W} < \sigma_{v}^{2} + t| \mathcal{H}_{1} \}
\\\leq&
\mathbb{P}\{ |E_{W} - \sigma_{v}^{2} - E_{U}| \geq E_{U} - t| \mathcal{H}_{1} \}
\\\leq&
\frac{4 E_{U} \sigma_{v}^{2} + 2 \sigma_{v}^{4}}{K (E_{U} - t)^{2}}
\\=&
\frac{4 E_{U} \sigma_{v}^{2} + 2 \sigma_{v}^{4}}{ \left(\sqrt{K} E_{U} - \sqrt{\frac{4}{\delta}} \sigma_{v}^{2} \right)^{2}}.
\end{split}
\end{align*}

%
Thus, to obtain $1-\delta$-detection for a given $\delta$, 
Willie sets $t = \frac{2 \sigma_{v}^{2}}{\sqrt{\delta K}}$, to achieve $\alpha \leq \delta/2$. Willie sets $K > K_{0}$, to achieve $\beta \leq \delta/2$, where,
\[
\begin{split}
K_{0} 
=&
\left(
\frac{2 \sqrt{2 E_{U} \sigma_{v}^{2} + \sigma_{v}^{4}} + 2 \sigma_{v}^{2}}{\sqrt{\delta} E_{U}}
\right)^{2}
\\=&
\frac{4}{\delta \cdot \text{SNR}^{2}} \left( 1 + \sqrt{1 + 2 \cdot \text{SNR}}\right)^{2},
\end{split}
\]

where $\text{SNR} \triangleq \frac{E_{U}}{\sigma_{v}^{2}}$. As a consequence, $\alpha + \beta \leq \delta$. 
\end{proof}

\subsection{Proof of \cref{theorem:An_inherently_stable_system_cant_be_stabilized_covertly_one_bit}}
\begin{proof}(\cref{theorem:An_inherently_stable_system_cant_be_stabilized_covertly_one_bit})
Willie is observing to Alice's system's output through a clean channel, when $|a| < 1$ and in steady state. Therefore, Willie's observation at time $n$, is an AR(1) signal controlled or not.
We prove the impossibility of covert control in this case, using the following detection method: 
Willie observes the system's output through a clean channel for $K+1$ samples, in any time sample $n$ Willie evaluates: $Y_{n} \triangleq X_{n} - aX_{n-1}$, and then calculates the average energy of $Y_{n}$, i.e., 
$E_{W} = \frac{1}{K}\sum_{n}^{}{Y_{n}^{2}}$. 
Willie compares $E_{W}$ to some expected energy level, in order to decide if the system is being controlled or not.

Willie preforms hypotheses testing approach to decide if Alice is controlling the system or not, he uses the following hypotheses,
\begin{align*}
\begin{matrix}
\mathcal{H}_{0}: &
Y_{n} = Z_{n}, & 
n = k, \ldots, K + k - 1,
\\
\mathcal{H}_{1}: &
Y_{n} = Z_{n} - U_{n}, &
n = k, \ldots, K + k - 1.
\end{matrix}
\end{align*}

Under the null hypothesis, Willie observes an i.i.d.\ process, thus, the mean and the variance of $E_{W}$ under the assumption that $\mathcal{H}_{0}$ is true are,
\begin{align}\label{eq:One_bit:Output:Mean_H_0}
\begin{split}
\mathbb{E} \left[ E_{W}| \mathcal{H}_{0} \right]
=&
\mathbb{E} \left[ \frac{1}{K}\sum_{n}^{}{ Y_{n}^{2}} | \mathcal{H}_{0} \right]
\\=&
\frac{1}{K}\sum_{n}^{}{ \mathbb{E} \left[ Z_{n}^{2}\right]} 
\\=&
\sigma_{Z}^{2},
\end{split}
\end{align}

\begin{align}\label{eq:One_bit:Output:Variance_H_0}
\begin{split}
\var{}{}{ E_{W}| \mathcal{H}_{0} }
=&
\var{}{}{ \frac{1}{K}\sum_{n}^{}{ Y_{n}^{2}}| \mathcal{H}_{0}}
\\=&
\frac{1}{K^{2}} \sum_{n}^{}{ \var{}{}{ Z_{n}^{2}}}
\\=&
\frac{1}{K^{2}} \sum_{n}^{}{ (m_{Z}(4) - \sigma_{Z}^{4})}
\\=&
\frac{m_{Z}(4) - \sigma_{Z}^{4}}{K},
\end{split}
\end{align}

where $m_{Z}(4)$ is the fourth moment of $Z$. 
In a similar fashion, under the alternative hypothesis, Willie observes an i.i.d.\ process with the control signal, which are both independent at the same time samples, thus, the mean and the variance of $E_{W}$ under the assumption that $\mathcal{H}_{1}$ is true are,
\begin{align}\label{eq:One_bit:Output:Mean_H_1}
\begin{split}
\mathbb{E} \left[ E_{W}| \mathcal{H}_{1} \right]
=&
\mathbb{E} \left[ \frac{1}{K}\sum_{n}^{}{ Y_{n}^{2}} | \mathcal{H}_{1} \right]
\\=&
\frac{1}{K}\sum_{n}^{}{ \mathbb{E} \left[ (Z_{n} - U_{n})^{2}\right]} 
\\=&
\frac{1}{K}\sum_{n}^{}{ (\expc{}{}{Z_{n}^{2}} + \expc{}{}{U_{n}^{2}})} 
\\=&
\sigma_{Z}^{2} + E_{U},
\end{split}
\end{align}

where $Z_{n} \indep U_{n}$ and $U_{n}^{2} = E_{U}$. 
\begin{align}\label{eq:One_bit:Output:Second_moment_H_1}
\begin{split}
\expc{}{}{ E_{W}^{2}| \mathcal{H}_{1} }
=&
\frac{1}{K^{2}} 
\sum_{n,m}^{}{ \expc{}{}{(Z_{n} - U_{n})^{2} (Z_{m} - U_{m})^{2}}}
\\\stackrel{(a)}{=}&
\frac{1}{K^{2}} \left( 
\sigma_{Z}^{4} K (K-1) + m_{Z}(4) K \right.
\\&+ \left.
2 E_{U} \sigma_{Z}^{2} K^{2} 
+
4 E_{U} \sigma_{Z}^{2} K + E_{U}^{2} K^{2}
\right)
\\=&
\left( \sigma_{Z}^{2} + E_{U} \right)^{2} 
\\&+ 
\frac{m_{Z}(4) - \sigma_{Z}^{4} + 4E_{U} \sigma_{Z}^{2}}{K},
\end{split}
\end{align}


where (a) is due to, 
\begin{align*}
\mathbb{E} &\left[ (Z_{n} - U_{n})^{2} (Z_{m} - U_{m})^{2} \right]
\\=&
\expc{}{}{(Z_{n}^{2} - 2 Z_{n} U_{n} + U_{n}^{2}) (Z_{m}^{2} - 2 Z_{m} U_{m} + U_{m}^{2})}
\\=&
\mathbb{E}\left[ Z_{n}^{2} Z_{m}^{2} - 2 Z_{n}^{2} Z_{m} U_{m} + Z_{n}^{2} U_{m}^{2}
\right.
\\&
-2 Z_{m}^{2} Z_{n} U_{n} + 4 Z_{n} U_{n} Z_{m} U_{m} - 2 Z_{n} U_{n} U_{m}^{2}
\\&
\left.
+ Z_{m}^{2} U_{n}^{2} - 2 Z_{m} U_{m} U_{n}^{2} + U_{n}^{2} U_{m}^{2}
\right]
\\=&
m_{Z}(4) \delta(n-m) + \sigma_{Z}^{4} (1 - \delta(n-m)) + E_{U} \sigma_{Z}^{2}
\\&
+ 4 E_{U} \sigma_{Z}^{2} \delta(n-m) 
+ 4 \expc{}{}{Z_{n} U_{n} Z_{m} U_{m}} (1 - \delta(n-m))
\\&+
E_{U} \sigma_{Z}^{2} + E_{U}^{2}
\\\stackrel{(b)}{=}&
m_{Z}(4) \delta(n-m) + \sigma_{Z}^{4} (1 - \delta(n-m)) + 2 E_{U} \sigma_{Z}^{2} 
\\& 
+ 4 E_{U} \sigma_{Z}^{2} \delta(n-m) + E_{U}^{2},
\end{align*}

where (b) is since, 
\begin{align*}
\begin{split}
\mathbb{E} & \left[ Z_{n} U_{n} Z_{m} U_{m}  \right] (1 - \delta(n-m))
\\=&
1_{m < n} \expc{}{}{Z_{n}} \expc{}{}{U_{n} Z_{m} U_{m}}
+
1_{m > n} \expc{}{}{Z_{m}} \expc{}{}{Z_{n} U_{n} U_{m}}
\\=&
0,
\end{split}
\end{align*}

summing the expression above yields, 
\begin{align*}
\begin{split}
\sum_{n,m}^{}{ \expc{}{}{(Z_{n} - U_{n})^{2} (Z_{m} - U_{m})^{2}}}  
=&
m_{Z}(4) K 
+
\sigma_{Z}^{4} K (K - 1) 
\\&+
2 E_{U} \sigma_{Z}^{2} K^{2} 
+
4 E_{U} \sigma_{Z}^{2} K 
\\&+
E_{U}^{2} K^{2}.
\end{split}
\end{align*}

Using \cref{eq:One_bit:Output:Mean_H_1,eq:One_bit:Output:Second_moment_H_1}, 
\begin{align}\label{eq:One_bit:Output:Variance_H_1}
\begin{split}
\var{}{}{ E_{W}| \mathcal{H}_{1} }
=&
\expc{}{}{ E_{W}^{2}| \mathcal{H}_{1} } - \expc{}{2}{ E_{W}| \mathcal{H}_{1} }
\\\stackrel{}{=}&
\frac{m_{Z}(4) - \sigma_{Z}^{4} + 4E_{U} \sigma_{Z}^{2}}{K},
\end{split}
\end{align}

If $\mathcal{H}_{0}$ is true, then $E_{W}$ should be close to $\mathbb{E}\left[ E_{W}| \mathcal{H}_{0} \right]$. Willie picks a threshold which we denote as $t$, and compares the value of $E_{W}$ to $\sigma_{Z}^{2} + t$. Willie accepts $\mathcal{H}_{0}$ if  $E_{W} < \sigma_{Z}^{2} + t$ and rejects it otherwise. 
We bound the false alarm probability using 
\cref{eq:One_bit:Output:Variance_H_0,eq:One_bit:Output:Mean_H_0} and with Chebyshev’s inequality,
\begin{align*}
\begin{split}
\alpha
=&
\mathbb{P}\{ E_{W} \geq \sigma_{Z}^{2} + t| \mathcal{H}_{0} \}
\\\leq&
\mathbb{P}\{ |E_{W} - \sigma_{Z}^{2}| \geq t| \mathcal{H}_{0} \}
\\\leq&
\frac{m_{Z}(4) - \sigma_{Z}^{4}}{K t^{2}}.
\end{split}
\end{align*}

Thus, to obtain $\alpha \leq \frac{\delta}{2}$, Willie sets $t = \sqrt{\frac{m_{Z}(4) - \sigma_{Z}^{4}}{K \delta / 2}}$.
The probability of a miss detection, $\beta$, is the probability that $E_{W} < \sigma_{Z}^{2} + t$ when $\mathcal{H}_{1}$ is true. We bound $\beta$ using \cref{eq:One_bit:Output:Variance_H_1,eq:One_bit:Output:Mean_H_1} and with Chebyshev’s inequality,
\begin{align*}
\begin{split}
\beta
=&
\mathbb{P}\{ E_{W} < \sigma_{Z}^{2} + t| \mathcal{H}_{1} \}
\\\leq&
\mathbb{P}\{ |E_{W} - (\sigma_{Z}^{2} + E_{U})| > E_{U} - t| \mathcal{H}_{1} \}
\\\leq&
\frac{m_{Z}(4) - \sigma_{Z}^{4} + 4E_{U} \sigma_{Z}^{2}}{K (E_{U} - t)^{2}}
\\=&
\frac{m_{Z}(4) - \sigma_{Z}^{4} + 4E_{U} \sigma_{Z}^{2}}{\left(\sqrt{K} E_{U} - \sqrt{\frac{m_{Z}(4) - \sigma_{Z}^{4}}{\delta / 2}} \right)^{2}}.
\end{split}
\end{align*}

Thus, to obtain $\beta \leq \frac{\delta}{2}$, Willie sets his observation window to be at least, 
\[
K_{0}
=
\frac{1}{E_{U}^{2}} \left( \sqrt{\frac{m_{Z}(4) - \sigma_{Z}^{4} + 4E_{U} \sigma_{Z}^{2}}{\delta/2}}
+
\sqrt{\frac{m_{Z}(4) - \sigma_{Z}^{4}}{\delta/2}}
\right)^{2},
\]

by doing so, Willie can detect with arbitrarily low error probability Alice's control actions with the one-bit controller, i.e., Willie achieves $1-\delta$-detection for any $\delta > 0$.
\end{proof}

\subsection{Proof of \cref{theorem:Achievable_gain_using_threshold_controller}}
In this subsection, the supporting claims of \cref{theorem:Achievable_gain_using_threshold_controller} and their proofs will be presented.

In order to prove \cref{theorem:Achievable_gain_using_threshold_controller}, we will use \cref{lemma:Connection_between_error_probabilities_to_total_variation,lemma:Connection_between_divergence_to_total_variation} to bound the sum of error probabilities. For doing so, we first give several supporting claims to upper bound the relevant Kullback-Leibler divergence.

Consider a specific case of a Gaussian AR(1) system, with and without the threshold controller. In this case, 
given $X_{k-1} = x_{k-1}$, $X_{k} \sim \mathcal{N}(a x_{k-1}, \sigma_{Z}^{2}), \quad \forall 2 \leq k \leq \tau$, where $\tau$ is the first threshold crossing time and $X_{0} = 0$. 

Denote $\textbf{X}^{(n)} \triangleq [X_{1}, X_{2}, \ldots, X_{n}]^{T}$. 
For $n \leq \tau$, i.e., without any control action, the probability density function of an AR(1) process, can be easily evaluated (see \cref{lemma:PDF_of_n_samples_AR_1_process}). 
First, we note the following. 
\begin{claim}\label{claim:Independency_of_AR_1_vectors}
Assume that an AR(1) system operates with the threshold controller, for a total of $N$ time samples, and there is only one crossing time which we indicate as $\tau_{1} < N$. 
Then, $\textbf{X}^{(1,\tau_{1})} |_{\tau = \tau_{1}} = [X_{1}, \ldots, X_{\tau_{1}}]^{T}$ and $\textbf{X}^{(\tau_{1} + 1,N)} |_{\tau = \tau_{1}} = [X_{\tau_{1} + 1}, \ldots, X_{N}]^{T}$ are two independent random vectors. 
\end{claim}

\begin{proof}
Since the threshold controller operates only at $\tau_{1} + 1$, which we assumed to be known, then $U_{\tau_{1} + 1} = a X_{\tau_{1}}$ and otherwise $U_{n} = 0$. 
I.e., the system operates undisturbed for $\tau_{1} < N$ time samples, which constitute an AR(1) process at that interval. On the other hand, after a correction has been made, the system state $X_{\tau_{1} + 1} = a X_{\tau_{1}} + Z_{\tau_{1} + 1} - U_{\tau_{1} + 1} = Z_{\tau_{1} + 1}$, have the same distribution as $X_{1} = Z_{1}$, since $Z_{n}$ is an i.i.d.\ process. 
Therefore, $\textbf{X}^{(1,\tau_{1})} |_{\tau = \tau_{1}}$ and $\textbf{X}^{(\tau_{1} + 1,N)} |_{\tau = \tau_{1}}$ have the same distribution, but with possibly different dimensions, and can be written as, \begin{align*}
\begin{split}
\textbf{X}^{(1,\tau_{1})} |_{\tau = \tau_{1}}
=&
\begin{pmatrix}
X_{1} \\
\vdots \\
X_{\tau_{1}}
\end{pmatrix}
\\=&
\begin{pmatrix}
1 & 0 & \cdots & 0 \\
\vdots & \vdots & \ddots & \vdots \\
a^{\tau_{1} - 1} & a^{\tau_{1} - 2} & \cdots & 1 
\end{pmatrix}
\begin{pmatrix}
Z_{1} \\
\vdots \\
Z_{\tau_{1}} \\
\end{pmatrix},
\\
\textbf{X}^{(\tau_{1} + 1,N)} |_{\tau = \tau_{1}}
=&
\begin{pmatrix}
X_{\tau_{1} + 1} \\
\vdots \\
X_{N}
\end{pmatrix}
\\=&
\begin{pmatrix}
1 & 0 & \cdots & 0 \\
\vdots & \vdots & \ddots & \vdots \\
a^{N - \tau_{1} - 1} & a^{N - \tau_{1} - 2} & \cdots & 1 
\end{pmatrix}
\begin{pmatrix}
Z_{\tau_{1} + 1} \\
\vdots \\
Z_{N}
\end{pmatrix},
\end{split}
\end{align*}

hence, $\textbf{X}^{(N)} |_{\tau = \tau_{1}}$ can be written as, 
\begin{align*}
\begin{split}
\textbf{X}^{(N)} |_{\tau = \tau_{1}}
=&
\begin{pmatrix}
\textbf{X}^{(1,\tau_{1})} \\
\textbf{X}^{(\tau_{1} + 1,N)}
\end{pmatrix}
\\=&
\begin{pmatrix}
\textbf{A}^{(\tau_{1})} & \textbf{0} \\
\textbf{0} & \textbf{A}^{(N-\tau_{1})}
\end{pmatrix}
\begin{pmatrix}
\textbf{Z}^{(1,\tau_{1})} \\
\textbf{Z}^{(\tau_{1} + 1,N)}
\end{pmatrix}.
\end{split}
\end{align*}

Now, $\textbf{Z}$ is an i.i.d.\ random vector, thus $\textbf{X}^{(1,\tau_{1})} \indep \textbf{X}^{(\tau_{1} + 1,N)}$.
\end{proof}

In terms of the density function, we have the following. 
\begin{claim}\label{claim:Cond_PDF_of_controlled_AR_1_process}
Consider a Gaussian AR(1) system. When using the threshold controller, the PDF of the system state vector $\textbf{X}^{(N)}$ conditioned to $\tau$, when there is only one crossing time, is given by
\begin{align}\label{eq:Cond_PDF_of_controlled_AR_1_process}
\begin{split}
f_{\textbf{X}^{(N)} |\tau = \tau_{1}}(\textbf{x}_{N} | \tau = \tau_{1})
=
\frac{1}{\left( 2 \pi \right)^{\frac{N}{2}} \sqrt{|\bs{\Sigma}_{N}|}} e^{-\frac{1}{2} \textbf{x}_{N}^{T} \bs{\Sigma}_{N}^{-1} \textbf{x}_{N}},
\end{split}
\end{align}

where $\tau_{1}$ is a real parameter, and $\bs{\Sigma}_{N}$ 
defined as follows, 
\begin{align*}
\begin{split}
\bs{\Sigma}_{N}
=&
\left(
\begin{matrix}
\bs{\Sigma}_{\tau_{1}} & \textbf{0} \\
\textbf{0} & \bs{\Sigma}_{N - \tau_{1}}
\end{matrix}
\right).
\end{split}
\end{align*}

I.e., $\textbf{X}^{(N)} \sim \mathcal{N}(\textbf{0}, \bs{\Sigma}_{N})$ and $\bs{\Sigma}_{N}$ is a covariance matrix of $N$ samples Gaussian AR(1) process.
\end{claim}

\begin{proof}
Since $\textbf{X}^{(1,\tau_{1})} \indep \textbf{X}^{(\tau_{1} + 1,N)}$, we have
\begin{align*}
\begin{split}
f_{\textbf{X}^{(N)} |\tau = \tau_{1}}(\textbf{x} | \tau = \tau_{1})
=&
f_{\textbf{X}^{(1,\tau_{1})} |\tau = \tau_{1}}(\textbf{x}_{1} | \tau = \tau_{1})
\\&
f_{\textbf{X}^{(N-\tau_{1},N)} |_{\tau = \tau_{1}}}(\textbf{x}_{2} | \tau = \tau_{1})
\\=&
\frac{1}{\left( 2 \pi \right)^{\frac{\tau_{1}}{2}} \sqrt{|\bs{\Sigma}_{\tau_{1}}|}} e^{-\frac{1}{2} \textbf{x}_{1}^{T} \bs{\Sigma}_{\tau_{1}}^{-1} \textbf{x}_{1}}
\\&
\frac{1}{\left( 2 \pi \right)^{\frac{N-\tau_{1}}{2}} \sqrt{|\bs{\Sigma}_{N-\tau_{1}}|}} e^{-\frac{1}{2} \textbf{x}_{2}^{T} \bs{\Sigma}_{N-\tau_{1}}^{-1} \textbf{x}_{2}}
\\=&
\frac{1}{\left( 2 \pi \right)^{\frac{N}{2}} \sqrt{|\bs{\Sigma}_{\tau_{1}}| |\bs{\Sigma}_{N-\tau_{1}}|}} 
\\&
e^{-\frac{1}{2} \left( \textbf{x}_{1}^{T} \bs{\Sigma}_{\tau_{1}}^{-1} \textbf{x}_{1} + \textbf{x}_{2}^{T} \bs{\Sigma}_{N-\tau_{1}}^{-1} \textbf{x}_{2}\right)}
\\\stackrel{(a)}{=}&
\frac{1}{\left( 2 \pi \right)^{\frac{N}{2}} \sqrt{|\bs{\Sigma}_{N}|}} e^{-\frac{1}{2} \textbf{x}_{N}^{T} \bs{\Sigma}_{N}^{-1} \textbf{x}_{N}},
\end{split}
\end{align*}

where (a) is by the relations: 
$\textbf{x}_{1} = \left[ x_{1}, \ldots, x_{\tau_{1}} \right]^{T}$, 
$\textbf{x}_{2} = \left[ x_{\tau_{1} + 1}, \ldots, x_{N} \right]^{T}$, $\textbf{x}_{N} = \left[ \textbf{x}_{1}^{T}, \textbf{x}_{2}^{T} \right]^{T}$, and 
\begin{align*}
\begin{split}
\bs{\Sigma}_{N}
=&
\left(
\begin{matrix}
\bs{\Sigma}_{\tau_{1}} & \textbf{0} \\
\textbf{0} & \bs{\Sigma}_{N - \tau_{1}}
\end{matrix}
\right).
\end{split}
\end{align*}
\end{proof}

\begin{coroc}
In general, assume that the system operates with the threshold controller for a total of $N$ time samples. Let us indicate 
$\{ \tau_{k} \}_{k = 1}^{m}$ as the series of crossing times. 
Therefore, \cref{eq:Cond_PDF_of_controlled_AR_1_process} is true for any $\tau_{k} + 1 \leq n \leq \tau_{k+1} \quad \forall k \in \mathbb{Z}_{+}$, where $\tau_{0} = 0$. 
In addition, as a result of the independence of the system state between each time interval, which no control action is been made, 
thus $\left\lbrace \textbf{X}^{(\tau_{k} + 1, \tau_{k+1})} | \bs{\tau} \right\rbrace_{k=1}^{m}$ are independent random vectors. Therefore, \cref{eq:Cond_PDF_of_controlled_AR_1_process} can be generalized to,
\begin{align}\label{eq:Cond_PDF_of_controlled_AR_1_process_general}
\begin{split}
f_{\textbf{X}^{(N)} | \bs{\tau}}(\textbf{x} | \bs{\tau})
=
%
\frac{1}{\left( 2 \pi \right)^{\frac{N}{2}} \sqrt{|\bs{\Sigma}|}} e^{-\frac{1}{2} \textbf{x}^{T} \bs{\Sigma}^{-1} \textbf{x}},
\end{split}
\end{align}

where $\bs{\Sigma}$ defined as follows,
\begin{align*}
\begin{split}
\bs{\Sigma}_{N}
=&
\begin{pmatrix}
\bs{\Sigma}_{\tau_{1}} & \\
 &  & \bs{\Sigma}_{\tau_{2}} &  \\
 & & & \ddots  \\
 & & & & \bs{\Sigma}_{N - \sum_{k}{\tau_{k}}}
\end{pmatrix}.
\end{split}
\end{align*}

%
%
\end{coroc}

We can now turn to the main technical claim, which bounds the divergence between the two relevant measures.
\begin{claim}\label{claim:Divergence_of_controlled_and_uncontrolled_AR_1}
Let $\textbf{X}^{(n)}$ be the vector of the system state of an uncontrolled Gaussian AR(1) system. Denote by $\tilde{\textbf{X}}^{(n)}$ the vector of the system state under one control action. The relative entropy between the PDF of $\textbf{X}^{(n)}$ and the PDF of $\tilde{\textbf{X}}^{(n)}$, is bounded by
\begin{align}\label{eq:Divergence_of_controlled_and_uncontrolled_AR_1}
\begin{split}
& D_{KL} \left( 
f_{\textbf{X}^{(n)}}(\textbf{x}) || f_{\tilde{\textbf{X}}^{(n)}}(\textbf{x})
\right)
\\\leq&
\frac{1}{2}
\expc{p_{\tau}(\tau_{1})}{}{
\tr{
\bs{\Sigma}_{\tilde{\textbf{X}}^{(n)} | \tau_{1}}^{-1} \bs{\Sigma}_{\textbf{X}^{(n)}}} 
+ \log{
\frac{\left\vert \bs{\Sigma}_{\tilde{\textbf{X}}^{(n)} | \tau_{1}} \right\vert}
{\left\vert \bs{\Sigma}_{\textbf{X}^{(n)}} \right\vert}}} 
-\frac{n}{2},
\end{split}
\end{align}

where $\bs{\Sigma}_{\tilde{\textbf{X}}^{(n)} |\tau_{1}}$ and $\bs{\Sigma}_{\textbf{X}^{(n)}}$ are the covariance matrices of $\tilde{\textbf{X}}^{(n)}  |_{\tau  = \tau_{1}}$ and $\textbf{X}^{(n)}$, respectively. 
\end{claim}

\begin{proof}
we have,
\begin{align*}
\begin{split}
D_{KL} &\left( 
f_{\textbf{X}^{(n)}} || f_{\tilde{\textbf{X}}^{(n)}}
\right)
\\\triangleq&
\mathbb{E}_{f_{\textbf{X}^{(n)}}}\left[
\log{\frac{f_{\textbf{X}^{(n)}}}{f_{\tilde{\textbf{X}}^{(n)}}}} 
\right]
\\=&
\expc{f_{\textbf{X}^{(n)}}}{}{
\log{\frac{f_{\textbf{X}^{(n)}}}{\sum_{\tau_{1} = 1}^{n}{f_{\tilde{\textbf{X}}^{(n)} | \tau}(\textbf{x} | \tau  = \tau_{1})
p_{\tau}(\tau_{1})}}}}
\\=&
\expc{f_{\textbf{X}^{(n)}}}{}{
\log{f_{\textbf{X}^{(n)}}} - \log{\sum_{\tau_{1} = 1}^{n}{f_{\tilde{\textbf{X}}^{(n)} | \tau}(\textbf{x} | \tau  = \tau_{1})
p_{\tau}(\tau_{1})}}}
\\\stackrel{(a)}{=}&
-H(f_{\textbf{X}^{(n)}})
- \expc{f_{\textbf{X}^{(n)}}}{}{\log{ \expc{p_{\tau}(\tau_{1})}{}{f_{\tilde{\textbf{X}}^{(n)} | \tau}(\textbf{x} | \tau  = \tau_{1})}}}
\\\stackrel{(b)}{\leq}&
-H(f_{\textbf{X}^{(n)}})
- \expc{f_{\textbf{X}^{(n)}}}{}{\expc{p_{\tau}(\tau_{1})}{}{\log{ f_{\tilde{\textbf{X}}^{(n)} | \tau}(\textbf{x} | \tau  = \tau_{1})}}}
\\\stackrel{(c)}{=}&
-H(f_{\textbf{X}^{(n)}})
- \expc{p_{\tau}(\tau_{1})}{}{\expc{f_{\textbf{X}^{(n)}}}{}{\log{ f_{\tilde{\textbf{X}}^{(n)} | \tau}(\textbf{x} | \tau  = \tau_{1})}}}
\\\stackrel{(d)}{=}&
\expc{p_{\tau}(\tau_{1})}{}{H_{CE} \left( f_{\textbf{X}^{(n)}}(\textbf{x}), f_{\tilde{\textbf{X}}^{(n)} | \tau}(\textbf{x} | \tau  = \tau_{1}) \right)}
-H(f_{\textbf{X}^{(n)}})
\\\stackrel{(e)}{=}&
\expc{p_{\tau}(\tau_{1})}{}{
\dive{KL}{f_{\textbf{X}^{(n)}}(\textbf{x})} {f_{\tilde{\textbf{X}}^{(n)} | \tau}(\textbf{x} | \tau  = \tau_{1})}}
\\\stackrel{(f)}{=}&
\frac{1}{2}
\expc{p_{\tau}(\tau_{1})}{}{
\tr{
\bs{\Sigma}_{\tilde{\textbf{X}}^{(n)} | \tau_{1}}^{-1} \bs{\Sigma}_{\textbf{X}^{(n)}}} 
+ \log{
\frac{\left\vert \bs{\Sigma}_{\tilde{\textbf{X}}^{(n)} | \tau_{1}} \right\vert}
{\left\vert \bs{\Sigma}_{\textbf{X}^{(n)}} \right\vert}}} 
-\frac{n}{2},
%
%
\end{split}
\end{align*}

where (a) is by $H(P) = - \expc{P}{}{\log{P}}$ as the entropy, and by representing the sum as an expectation according to $p_{\tau}$. (b) is due to Jensen's inequality. 
(c) is due to changing the order of the expectations, since $\expc{p_{\tau}}{}{\cdot}$ is a discrete and finite expectation. 
(d) is by the definition of the \emph{cross-entropy}, 
$H_{CE}(P,Q) = -\expc{P}{}{\log{Q}}$. 
(e) is due to the following relation: $\dive{KL}{P}{Q} = H_{CE}(P,Q) - H(P)$, and since $H(f_{\textbf{X}^{(n)}})$ is independent of $\tau$. (f) is by \cref{lemma:Divergence_of_gaussian_vector}. 
%
\end{proof}





\begin{claim}\label{claim:Divergence_of_controlled_and_uncontrolled_AR_1_of_stable_system_in_steady_state}
Consider a Gaussian AR(1) system with $|a| < 1$, and the system is in steady state, then
\begin{align}\label{eq:Divergence_of_controlled_and_uncontrolled_AR_1_of_stable_system_in_steady_state}
\begin{split}
\dive{KL}{f_{\textbf{X}^{(n)}}(\textbf{x})} {f_{\tilde{\textbf{X}}^{(n)} | \tau}(\textbf{x} | \tau  = \tau_{1})}
=&
\frac{1}{2} \log{\frac{1}{1 - a^{2}}}.
\end{split}
\end{align}

\end{claim}

\begin{corollary}\label{corollary:Diveragence_upper_bound}
For a Gaussian AR(1) system with $|a| < 1$ and in steady state, by 
\cref{claim:Divergence_of_controlled_and_uncontrolled_AR_1_of_stable_system_in_steady_state,claim:Divergence_of_controlled_and_uncontrolled_AR_1},
we get, 
\[
\dive{KL}{f_{\textbf{X}^{(n)}}}{f_{\tilde{\textbf{X}}^{(n)}}}
\leq 
\frac{1}{2} \log
\frac{1}{1 - a^{2}}, 
\]

since 
$\dive{KL}{f_{\textbf{X}^{(n)}}(\textbf{x})} {f_{\tilde{\textbf{X}}^{(n)} | \tau}(\textbf{x} | \tau_{1})},$ 
does not depend on $\tau_{1}$.
\end{corollary}

\begin{proof}(\cref{claim:Divergence_of_controlled_and_uncontrolled_AR_1_of_stable_system_in_steady_state})
We need to show two things. First, we will prove that,
\[
\tr{\bs{\Sigma}_{\tilde{\textbf{X}}^{(n)} | \tau_{1}}^{-1} \bs{\Sigma}_{\textbf{X}^{(n)}}} = n,
\]

and second, 
\[
\log{\frac{\left\vert \bs{\Sigma}_{\tilde{\textbf{X}}^{(n)} | \tau_{1}} \right\vert}
{\left\vert \bs{\Sigma}_{\textbf{X}^{(n)}} \right\vert}} 
=
\log{\frac{1}{1 - a^{2}}}.
\]

Recall that an $n$-tupple of a Gaussian AR(1) process with $|a| < 1$ at steady state, has the following covariance matrix
$[\bs{\Sigma}]_{i,j} = \frac{\sigma_{Z}^{2}}{1 - a^{2}} a^{|i-j|}$. 
Hence, in our case, 
\[
\left[  \bs{\Sigma}_{\textbf{X}^{(n)}} \right]_{i,j}
= 
\frac{\sigma_{Z}^{2}}{1 - a^{2}} a^{|i-j|}, \quad 1 \leq i,j \leq n,
\]

and, 
\[
\bs{\Sigma}_{\tilde{\textbf{X}}^{(n)} | \tau_{1}}
=
\begin{pmatrix}
\bs{\Sigma}_{\tau_{1}} & \textbf{0}_{\tau_{1}, n - \tau_{1}} \\
\textbf{0}_{n - \tau_{1}, \tau_{1}} & \bs{\Sigma}_{n - \tau_{1}},
\end{pmatrix},
\]

where,
\[
\begin{split}
\left[ \bs{\Sigma}_{\tau_{1}} \right]_{i,j}
=&
\frac{\sigma_{Z}^{2}}{1 - a^{2}} a^{|i-j|}, \quad 1 \leq i,j \leq \tau_{1},
\\
\left[ \bs{\Sigma}_{n-\tau_{1}} \right]_{i,j}
=&
\frac{\sigma_{Z}^{2}}{1 - a^{2}} a^{|i-j|}, \quad 1 \leq i,j \leq n-\tau_{1}.
\end{split}
\]

One can divide 
$\bs{\Sigma}_{\textbf{X}^{(n)}}$ 
to blocks the same way as 
$\bs{\Sigma}_{\tilde{\textbf{X}}^{(n)} | \tau_{1}}$, 
thus, 
\[
\bs{\Sigma}_{\textbf{X}^{(n)}}
=
\begin{pmatrix}
\bs{\Sigma}_{11} & \bs{\Sigma}_{12} \\
\bs{\Sigma}_{12}^{T} & \bs{\Sigma}_{22}
\end{pmatrix},
\]

where, 
\[
\begin{split}
\left[ \bs{\Sigma}_{11} \right]_{i,j}
=&
\left[ \bs{\Sigma}_{\tau_{1}} \right]_{i,j},
\\
\left[ \bs{\Sigma}_{22} \right]_{i,j}
=&
\left[ \bs{\Sigma}_{n-\tau_{1}} \right]_{i,j}.
\end{split}
\]

Therefore,
\[
\begin{split}
\bs{\Sigma}_{\tilde{\textbf{X}}^{(n)} | \tau_{1}}^{-1} \bs{\Sigma}_{\textbf{X}^{(n)}}
=&
\begin{pmatrix}
\bs{\Sigma}_{\tau_{1}}^{-1} & \textbf{0}_{\tau_{1}, n - \tau_{1}} \\
\textbf{0}_{n - \tau_{1}, \tau_{1}} & \bs{\Sigma}_{n - \tau_{1}}^{-1}
\end{pmatrix}
\begin{pmatrix}
\bs{\Sigma}_{11} & \bs{\Sigma}_{12} \\
\bs{\Sigma}_{12}^{T} & \bs{\Sigma}_{22}
\end{pmatrix}
\\=&
\begin{pmatrix}
\bs{\Sigma}_{\tau_{1}}^{-1} \bs{\Sigma}_{11} & \bs{\Sigma}_{\tau_{1}}^{-1} \bs{\Sigma}_{12} \\
\bs{\Sigma}_{n-\tau_{1}}^{-1} \bs{\Sigma}_{12}^{T} & \bs{\Sigma}_{n-\tau_{1}}^{-1} \bs{\Sigma}_{22}
\end{pmatrix}
\\=&
\begin{pmatrix}
\textbf{I}_{\tau_{1}} & \bs{\Sigma}_{\tau_{1}}^{-1} \bs{\Sigma}_{12} \\
\bs{\Sigma}_{n-\tau_{1}}^{-1} \bs{\Sigma}_{12}^{T} & \textbf{I}_{n-\tau_{1}} 
\end{pmatrix},
\end{split}
\]

which yields,
\[
\tr{\bs{\Sigma}_{\tilde{\textbf{X}}^{(n)} | \tau_{1}}^{-1} \bs{\Sigma}_{\textbf{X}^{(n)}}} = n.
\]

Now, we move to the second part of the proof. 
We use a LU decomposition for the matrix defined as
$\left[ \hat{\textbf{A}}_{n} \right]_{i,j} = a^{|i-j|}, \quad 1 \leq i,j \leq n$. 
Where the lower triangular matrix $\textbf{L}_{n}$ defined as 
$\left[ \textbf{L}_{n} \right]_{i,j} = a^{i-j} u(i-j), \quad 1 \leq i,j \leq n$, 
and the upper triangular matrix $\textbf{U}_{n}$ defined as 
$\left[ \textbf{U}_{n} \right]_{i,j} = a^{j-i} (1 - a^{2})^{(1 - \delta(i-1))} u(j-i), \quad 1 \leq i,j \leq n$. 
Where $\delta(\cdot)$ denote the Kronecker's delta, and $u(\cdot)$ is the discrete step function. If so, 
\[
\begin{split}
\left[ \hat{\textbf{A}}_{n} \right]_{i,j}
=&
\sum_{k=1}^{n}{\left[ \textbf{L}_{n} \right]_{i,k} \left[ \textbf{U}_{n} \right]_{k,j}}
\\=&
\sum_{k=1}^{n}{a^{i-k} u(i-k) a^{j-k} (1 - a^{2})^{(1 - \delta(k-1))} u(j-k)}
\\\stackrel{(a)}{=}&
a^{i+j}
\sum_{k=1}^{\gamma}{a^{-2k} (1 - a^{2})^{(1 - \delta(k-1))}}
\\\stackrel{(b)}{=}&
a^{i+j}
\left(  
a^{-2} + 
(1 - a^{2}) \sum_{k=2}^{\gamma}{a^{-2k}}
\right)
\\\stackrel{(c)}{=}&
a^{i+j}
\left(  
a^{-2} + 
(1 - a^{2}) \frac{a^{-2 \gamma} - a^{-2}}{1 - a^{2}}
\right)
\\=&
a^{i+j}
\left(  
a^{-2} + a^{-2 \gamma} - a^{-2}
\right)
\\=&
a^{i+j -2 \gamma}
\\\stackrel{(d)}{=}&
a^{|i-j|},
\end{split}
\]

where (a) is by the fact that $k \leq i$ and $k \leq j$, hence $k \leq \min{(i,j)} \triangleq \gamma$. 
(b) is by dividing the sum for $k = 1$ and for $k > 1$. 
(c) is due to sum of geometric series. 
(d) is by $i + j - 2\min{(i,j)} = |i-j|$. 

Now, consider that $\left[ \textbf{L}_{n} \right]_{i,i} = 1$ and 
$\left[ \textbf{U}_{n} \right]_{i,i} = \delta(i-1) + (1 - a^{2}) u(i-2)$, 
which yields, 
\[
|\hat{\textbf{A}}_{n}|
=
|\textbf{L}_{n}| |\textbf{U}_{n}|
=
1 \cdot 1 \cdot (1 - a^{2})^{n-1}
= 
(1 - a^{2})^{n-1},
\]

therefore, for $1 \leq \tau_{1} < n$ we have, 
\[
\begin{split}
& \log
\frac{\left\vert \bs{\Sigma}_{\tilde{\textbf{X}}^{(n)} | \tau_{1}} \right\vert}
{\left\vert \bs{\Sigma}_{\textbf{X}^{(n)}} \right\vert} 
\\=&
\log
\frac{\left\vert \frac{\sigma_{Z}^{2}}{1 - a^{2}} \hat{\textbf{A}}_{\tau_{1}} \right\vert \left\vert \frac{\sigma_{Z}^{2}}{1 - a^{2}} \hat{\textbf{A}}_{n - \tau_{1}} \right\vert}
{\left\vert \frac{\sigma_{Z}^{2}}{1 - a^{2}} \hat{\textbf{A}}_{n} \right\vert} 
\\=&
\log
\frac{\left( \frac{\sigma_{Z}^{2}}{1 - a^{2}} \right)^{\tau_{1}}
(1 - a^{2})^{\tau_{1}-1}
\left( \frac{\sigma_{Z}^{2}}{1 - a^{2}} \right)^{n - \tau_{1}}
(1 - a^{2})^{n - \tau_{1}-1}}
{\left( \frac{\sigma_{Z}^{2}}{1 - a^{2}} \right)^{n}
(1 - a^{2})^{n-1}} 
\\=&
\log
\frac{1}{1 - a^{2}}.
\end{split}
\]

\end{proof}

%
%
%

We can now give the proof of  \cref{theorem:Achievable_gain_using_threshold_controller}, which asserts that if the gain, $a$, is small enough, then Willie which observes the system's output, cannot distinguish if any
control action is been made by Alice. I.e., Willie cannot know if the system operates for $n$
samples without interference, or that Alice performes a reset.

\begin{proof}(\cref{theorem:Achievable_gain_using_threshold_controller})
We know that, 
\[
\begin{split}
\alpha + \beta
\stackrel{(a)}{=}&
1 - 
\mathcal{V}_{T}\left( f_{\textbf{X}^{(n)}}, f_{\tilde{\textbf{X}}^{(n)}} \right)
\\ \stackrel{(b)}{\geq}&
1 - 
\sqrt{\frac{1}{2} \dive{KL}{f_{\textbf{X}^{(n)}}}{f_{\tilde{\textbf{X}}^{(n)}}}}
\\ \stackrel{(c)}{\geq}&
1 - 
\sqrt{\frac{1}{4} \log\frac{1}{1 - a^{2}}}
\\ \stackrel{(d)}{\geq}&
1 - 
\varepsilon,
\end{split}
\]

(a) 
is due to 
\cref{lemma:Connection_between_error_probabilities_to_total_variation} . (b) is due to 
\cref{lemma:Connection_between_divergence_to_total_variation}. 
(c) 
by \cref{claim:Divergence_of_controlled_and_uncontrolled_AR_1_of_stable_system_in_steady_state},  and (d) is by the covertness criterion. 
Now we have, 
\[
\varepsilon
\geq
\sqrt{\frac{1}{4} \log\frac{1}{1 - a^{2}}},
\]

which yields,
\[
|a| 
\leq 
\sqrt{1 - 2^{-4 \varepsilon^{2}}}.
\]
\end{proof}

\subsection{Proof of \cref{theorem:Converse_conditional_distribution}}
\begin{proof}


Consider the following hypotheses testing problem that Willie uses, in order to detect Alice's control actions using the threshold controller.
\begin{align}\label{eq:Converse_conditional_distribution_hypotheses}
\begin{split}
\begin{matrix}
\mathcal{H}_{0}: X_{\tau+1} & \sim & \mathcal{N}\left( 0, \frac{\sigma_{Z}^{2}}{1 - a^{2}} \right)
\\
\mathcal{H}_{1}: X_{\tau+1} & \sim & \mathcal{N}(0, \sigma_{Z}^{2})
\end{matrix}
, \quad 1 \leq \tau \leq N.
\end{split}
\end{align}

Under $\mathcal{H}_{0}$, i.e., the system was not interrupted, the distribution of $X_{\tau+1}$ is given by the distribution of a steady state  AR(1) process, which is $\mathcal{N}\left( 0, \frac{\sigma_{Z}^{2}}{1 - a^{2}} \right)$ (by \cref{lemma:PDF_of_n_samples_AR_1_process}). 
Under $\mathcal{H}_{1}$, i.e., the system was interrupted, the distribution of $X_{\tau+1}$ is given by the distribution of the noise at $\tau+1$, which is $Z_{\tau+1} \sim \mathcal{N}\left( 0, \sigma_{Z}^{2} \right)$. The log-likelihood ratio is thus given by
\[
\begin{split}
\log{T'}
= & 
\log{\frac{f_{X_{\tau + 1}}(x_{\tau + 1} | \mathcal{H}_{1})}{f_{X_{\tau + 1}}(x_{\tau + 1} | \mathcal{H}_{0})}}
\\=&
\log{\frac{\sigma_{0}}{\sigma_{1}}}
-
\frac{x_{\tau + 1}^{2}}{2}\left( \frac{1}{\sigma_{1}^{2}} - \frac{1}{\sigma_{0}^{2}}\right)
\\=&
\log{\frac{1}{\sqrt{1 - a^{2}}}}
-
\frac{a^{2}}{2 \sigma_{Z}^{2}} x_{\tau + 1}^{2}.
\end{split}
\]

Hence, the equivalent test statistics is $T = X_{\tau+1}^{2}$, and the resulting decision rule is to compare $X_{\tau+1}^{2}$ to a threshold $t$, in particular 
\begin{align*}
T
=
x_{\tau + 1}^{2}
\DR{\mathcal{H}_{0}}{\mathcal{H}_{1}} & 
\frac{2 \sigma_{Z}^{2}}{a^{2}} \left( \log{\frac{1}{\sqrt{1 - a^{2}}}} - \log{t'} \right)
=
t.
\end{align*}

Since $X_{\tau+1} / \sigma_{k} \sim \mathcal{N}(0,1)$ for $k = 0,1$, hence, $(X_{\tau+1} / \sigma_{k})^{2} \sim \chi_{1}^{2}$, the false alarm probability is
\[
\begin{split}
\alpha
\triangleq &
\prob{X_{\tau+1}^{2} \leq t | \mathcal{H}_{0}}
\\=&
\prob{\frac{X_{\tau+1}^{2}}{\sigma_{0}^{2}} \leq \frac{t}{\sigma_{0}^{2}} | \mathcal{H}_{0}}
\\\stackrel{(a)}{=}&
1 - Q_{\chi_{1}^{2}}\left( \frac{t (1 - a^{2})}{\sigma_{Z}^{2}} \right)
\\\stackrel{(b)}{=}&
1 - 2Q\left( \frac{\sqrt{t (1 - a^{2})}}{\sigma_{Z}} \right),
\end{split}
\]

where (a) is since $Q_{\chi_{1}^{2}}(\cdot)$ is the right tail probability of $\chi_{1}^{2}$.  
(b) is by the relation $Q_{\chi_{1}^{2}}(z) = 2Q(\sqrt{z})$ \cite{Kay_detection}.
Applying the detection constraint, i.e., bounding $\alpha \leq \frac{\delta}{2}$, the threshold can be chosen to be
\[
\begin{split}
t
=
\frac{\sigma_{Z}^{2}}{1 - a^{2}}
\left( Q^{-1}\left( \frac{1 - \delta/2}{2} \right) \right)^{2}.
\end{split}
\]

On the other hand, the miss detection probability is
\[
\begin{split}
\beta
\triangleq &
\prob{X_{\tau+1}^{2} > t | \mathcal{H}_{1}}
\\=&
\prob{\frac{X_{\tau+1}^{2}}{\sigma_{1}^{2}} > \frac{t}{\sigma_{1}^{2}} | \mathcal{H}_{1}}
\\\stackrel{}{=}&
Q_{\chi_{1}^{2}}\left( \frac{t}{\sigma_{Z}^{2}} \right)
\\\stackrel{}{=}&
2Q\left( \frac{\sqrt{t}}{\sigma_{Z}} \right)
\\\stackrel{(a)}{=}&
2Q\left( \frac{1}{\sqrt{1 - a^{2}}} Q^{-1}\left( \frac{1 - \delta/2}{2} \right) \right)
\\\stackrel{(b)}{\leq}&
e^{-\frac{1}{2} \frac{1}{1 - a^{2}} \left(Q^{-1}\left( \frac{1 - \delta/2}{2} \right) \right)^{2}},
\end{split}
\]

where (a) is by substituting $t$. (b) is by $Q(x) \leq \frac{1}{2} e^{-\frac{x^{2}}{2}}$ \cite{Q_func}, 
and the fact that the argument is positive. By bounding $\beta \leq \frac{\delta}{2}$, we have
\[
|a|
\geq
\sqrt{1 - \frac{\left( Q^{-1}\left( \frac{1 -  \delta/2}{2} \right) \right)^{2}}{2 \log{\frac{2}{\delta}}}}.
\]
\end{proof}

\end{document}